\documentclass{article}

\usepackage{varioref}
\usepackage{srcltx}
\input{xy}
\xyoption{all}

\usepackage{varioref}
\usepackage{srcltx}
\input{xy}
\xyoption{all}

\usepackage{verbatim}
\usepackage{varioref}
\usepackage{amsfonts,
enumerate,
amsthm,
amsmath,latexsym,amscd,amssymb,doc}
\usepackage[latin1]{inputenc}
\usepackage{graphicx}
\usepackage{hyperref}
	
\usepackage{cyrillic}
\newcommand{\cyrrm}{\fontencoding{OT2}\selectfont\textcyrup}

\DeclareMathOperator{\Beilinson}{{\cyrrm{B}}} %for the Beilinson spectrum

		\newcommand{\category}[1]{\mathbf{#1}}
	
		\newcommand{\Com}{\category{Com}} % complexe
    \newcommand{\D}{\category {D}} % abgeleitete categ

    \newcommand{\Ab}{\category {Ab}} % abelsche Gruppen
    \newcommand{\Sets}{\category {Sets}} % Mengen
    \newcommand{\Top}{\category {Top}} % topologische Räume
    \newcommand{\Sm}{\category{Sm}}
    \newcommand{\AnSm}{\category{AnSm}} % analytic smooth spaces
		\newcommand{\Spt}{\category{Spt}} % spectra
		\newcommand{\PSh}{\category{PSh}} % presheaves

		\newcommand{\Shvv}{\category{Shv}} % sheaves
    \newcommand{\DM}{\category {DM}}
    
    \newcommand{\DMBei}{\DM_{\Beilinson}}

\newcommand{\Mod}{\category {Mod}} % modules

\newcommand{\SH}{\category{SH}} % stable homotopy category
\newcommand{\HoC}{\category{Ho}} % homotopy category (nonstable)
\newcommand{\HoCsect}{\HoC_{\mathrm{sect}, \bullet}} % homotopy category (nonstable)
	    \newtheoremstyle{Normal}% name
  {}%9pt}%      Space above, empty = `usual value'
  {}%9pt}%      Space below
  {}%\itshape}% Body font
  {}%         Indent amount (empty = no indent, \parindent = para indent)
  {\bfseries}% Thm head font
  {.}%        Punctuation after thm head
  { }%\newline}% Space after thm head: \newline = linebreak
  {}%         Thm head spec

    \theoremstyle{Normal} % Normal

    \newtheorem{Defi}{Definition}[section]

    \newtheorem{Conj}[Defi]{Conjecture}
    \newtheorem{Bem}[Defi]{Remark}
    \newtheorem{Bsp}[Defi]{Example}
     
    \newtheorem{Axio}[Defi]{Axiom}
    \newtheorem{Ques}[Defi]{Question}
    \newtheorem{Comp}[Defi]{Complements}
    
		\theoremstyle{remark}
	
    \theoremstyle{plain} %kursiv

    \newtheorem{Satz}[Defi]{Proposition}
    \newtheorem{DefiTheo}[Defi]{Definition and Theorem}
    \newtheorem{Theo}[Defi]{Theorem}
    \newtheorem{Folg}[Defi]{Corollary}
    \newtheorem{Lemm}[Defi]{Lemma}
    \newtheorem{DefiLemm}[Defi]{Definition and Lemma}

\newcommand{\refit}[1]{(\ref{item_#1})}

\newcommand{\refsect}[1]{Section \ref{sect_#1}}
\newcommand{\conj}{\begin{Conj}} 			\newcommand{\xconj}{\end{Conj}}											    
\newcommand{\ques}{\begin{Ques}} 			\newcommand{\xques}{\end{Ques}}											    
\newcommand{\axio}{\begin{Axio}} 			\newcommand{\xaxio}{\end{Axio}}											    
\newcommand{\bem}{\begin{Bem}} 			\newcommand{\xbem}{\end{Bem}}											    
\newcommand{\rema}{\begin{Bem}} 			\newcommand{\xrema}{\end{Bem}}											    \newcommand{\refre}[1]{Remark \ref{rema_#1}}
\newcommand{\defi}{\begin{Defi}} 			\newcommand{\xdefi}{\end{Defi}}										\newcommand{\refde}[1]{Definition \ref{defi_#1}}
\newcommand{\defitheo}{\begin{DefiTheo}} \newcommand{\xdefitheo}{\end{DefiTheo}} 
\newcommand{\defilemm}{\begin{DefiLemm}} \newcommand{\xdefilemm}{\end{DefiLemm}} \newcommand{\refdele}[1]{Definition and Lemma \ref{defilemm_#1}}

\newcommand{\lemm}{\begin{Lemm}}			\newcommand{\xlemm}{\end{Lemm}}											\newcommand{\refle}[1]{Lemma \ref{lemm_#1}}
\newcommand{\comp}{\begin{Comp}}			\newcommand{\xcomp}{\end{Comp}}											

\newcommand{\satz}{\begin{Satz}}			\newcommand{\xsatz}{\end{Satz}}										

\newcommand{\prop}{\begin{Satz}}			\newcommand{\xprop}{\end{Satz}}										
\newcommand{\theo}{\begin{Theo}}			\newcommand{\xtheo}{\end{Theo}}											\newcommand{\refth}[1]{Theorem \ref{theo_#1}}
\newcommand{\bsp}{\begin{Bsp}}				\newcommand{\xbsp}{\end{Bsp}}												
\newcommand{\exam}{\begin{Bsp}}				\newcommand{\xexam}{\end{Bsp}}												
\newcommand{\folg}{\begin{Folg}}				\newcommand{\xfolg}{\end{Folg}}
\newcommand{\coro}{\folg}				\newcommand{\xcoro}{\xfolg} 			\newcommand{\refcor}[1]{Corollary \ref{coro_#1}}
\newcommand{\cor}{\begin{Folg}}				\newcommand{\xcor}{\end{Folg}}
\newcommand{\mycomment}{\begin{comment}}				\newcommand{\xcomment}{\end{comment}}

\newcommand{\eqnarra}{\begin{eqnarray}}				\newcommand{\xeqnarra}{\end{eqnarray}}
\newcommand{\eqnarr}{\begin{eqnarray*}}				\newcommand{\xeqnarr}{\end{eqnarray*}}

\newcommand{\eqn}{\begin{equation}} 		\newcommand{\xeqn}{\end{equation}}
\newcommand{\refeq}[1]{(\ref{eqn_#1})}

\newcommand{\mylabel}[1]{\label{#1}}%\framebox[1.1\width]{$#1$}} %%% IF YOU GET BOTHERED by the frames, just comment out (but please don't delete) the \framebox part
\newcommand{\todo}[1]{\framebox[1.1\width]{#1}}
\newcommand{\status}[1]{}%\framebox[1.1\width]{STATUS $#1$}}

\newcommand{\explainiso}[1]{\stackrel{\text{#1}}{\cong}}

\newcommand{\Ho}{\category{Ho}}
\newcommand{\OF}{{\mathcal{O}_F}}

\newcommand{\dual}{^\vee}%{\check{\hspace{0.5cm}}}}%{\mathrm{dual}}} % dual vector space etc.
\newcommand{\twi}[1]{\{#1\}}

\renewcommand{\log}{\mathrm{log}\,}

\newcommand{\R}{\mathrm{R}} % right derived functor
\newcommand{\Hw}{\H_\mathrm{w}} %  weak Hodge cohomology
\newcommand{\HD}{\H_\mathrm{D}} %  Deligne cohomology
\newcommand{\HBetti}{\H_\mathrm{B}} %  Betti cohomology

\newcommand{\lr}{{\longrightarrow}}
\renewcommand{\r}{\rightarrow}

\renewcommand{\t}{{\otimes}}

\newcommand{\A}[1][1]{\mathbb{A}^{#1}}
\renewcommand{\P}[1][1]{\mathbb P^{#1}}
\newcommand{\N}{\mathbb{N}}
\newcommand{\Q}{\mathbb{Q}}

\newcommand{\An}{{\operatorname{An}}} % analytification

\newcommand{\CC}{\mathbb{C}}
\newcommand{\RR}{\mathbb{R}}
\newcommand{\Z}{\mathbb{Z}}

\newcommand{\Fp}{\mathbb{F}_p}

\renewcommand{\H}{\mathrm{H}}

\newcommand{\x}{{\times}}

\newcommand{\one}{\mathbf{1}}
\newcommand{\onehat}{\widehat \one{}}

\newcommand{\Beweis}{{\normalfont} \textbf{Proof}}

\newcommand{\lc}{\textit{loc.~cit.}}

\newcommand{\lcs}{\textit{loc.~cit.}\ }
\newcommand{\ocs}{\textit{op.~cit.}\ }
\newcommand{\Spec}{\mathrm{Spec} \; \! }
\newcommand{\SpecOF}{{\Spec} {\OF}}

\newcommand{\SpecZ}{{\Spec} {\Z}}

\newcommand{\SpecFp}{{\Spec} {\Fp}}

\newcommand{\Fr}{\operatorname{Fr}} % Frobenius
\newcommand{\Td}{\operatorname{Td}} % Todd class
\newcommand{\op}{\mathrm{op}}
\newcommand{\Gm}{\mathbb{G}_m}

\newcommand{\gr}{\operatorname{gr}} % grading

\newcommand{\tr}{\mathrm{tr}}
\newcommand{\ch}{\operatorname{ch}} % Chern character etc.
\newcommand{\cl}{\operatorname{cl}} % Chern character etc.
\newcommand{\Hom}{\mathrm{Hom}}
\newcommand{\End}{\mathrm{End}}
\newcommand{\IHom}{\underline{\Hom}}
\newcommand{\id}{\mathrm{id}}
\renewcommand{\gr}{\operatorname{gr}} % graduation

\newcommand{\im}{\operatorname{im}} %image
\newcommand{\M}{{\operatorname{M}}} % motive
\newcommand{\CH}{\mathrm{CH}}
\newcommand{\CHhat}{\widehat \CH{}}
\newcommand{\CHbar}{\overline \CH{}}
\newcommand{\Khat}{\widehat{\mathrm{K}}}
\newcommand{\Kbar}{\overline {K}}
\newcommand{\Shat}{\widehat {S}}
\newcommand{\Ghat}{\widehat {G}}
\newcommand{\Phat}{\widehat {P}}
\newcommand{\CHhatGS}{\widehat{\CH}_{\mathrm{GS}{}}}

\newcommand{\Gr}{\operatorname{Gr}} % Grassmannian

\newcommand{\cone}{\operatorname{cone}}

\newcommand{\Id}{\mathrm{Id}}
\newcommand{\colim}{{\operatornamewithlimits{colim}}}

\newcommand{\loc}[2]{[}

\newcommand{\pr}{\begin{proof}[\Beweis: ]}
\newcommand{\pf}{\pr}
\newcommand{\xpf}{\end{proof}}

\theoremstyle{definition}

\theoremstyle{definition}

\theoremstyle{definition}

\theoremstyle{definition}

\theoremstyle{definition}

\theoremstyle{definition}

\theoremstyle{definition}

\theoremstyle{definition}

\theoremstyle{definition}

\theoremstyle{definition}

\theoremstyle{definition}

\theoremstyle{definition}

\theoremstyle{definition}

\theoremstyle{definition}

\newcommand{\HB}{\mathrm{H}_{\Beilinson}}
\newcommand{\HBRR}{\mathrm{H}_{\Beilinson, \RR}}

\newcommand{\HBhat}{\widehat{\HB}}
\newcommand{\HBbar}{\overline{\HB}}
\newcommand{\HBhatRR}{\widehat{\HBRR}}
\newcommand{\HBhatS}{\widehat{\H_{{\Beilinson}, S}}}

\newcommand{\HBS}{\mathrm{H}_{{\Beilinson}, S }}

\newcommand{\Hhat}{\widehat{\mathrm{H}}}
\newcommand{\Hbar}{\overline{\mathrm{H}}}

\newcommand{\hofib}{\operatorname{hofib}} % homotpy fiber
\newcommand{\holim}{\operatorname{holim}} % homotopy limit

\newcommand{\cB}{\mathcal{B}}

\newcommand{\cD}{\mathcal{D}}
\newcommand{\cK}{\mathcal{K}}
\newcommand{\cO}{\mathcal{O}}

\newcommand{\cDz}{\H_\Deligne}

\newcommand{\Deligne}{\mathrm{D}}

\newcommand{\BGL}{\mathrm{BGL}}
\newcommand{\BGLhat}{\widehat{\BGL}}
\newcommand{\BGLbar}{\overline{\BGL}}

\newcommand{\SHnaive}{\SH^\mathrm{naive}}
\newcommand{\BGLnaive}{\BGL^\mathrm{naive}}

\usepackage{titlesec}
\usepackage[]{titletoc}

\titlecontents{subsubsubsection}[9em]{}{\contentslabel{3.9em}}%
{\hspace*{-1.2em}}{\titlerule*[0.675pc]{.}\contentspage}

\makeatletter
\newcounter{subsubsubsection}[subsubsection]
\renewcommand{\thesubsubsubsection}{\thesubsubsection.\@arabic\c@subsubsubsection}

\titleclass{\subsubsubsection}{straight}[\subsubsection]
\titleformat{\subsubsubsection}{\bf}{\thetitle}{1em}{}[]						
\titlespacing{\subsubsubsection}{0pt}{3.25ex plus 1ex minus 0.2ex}{1.5ex plus 0.2ex}

\makeatother

\usepackage{xr}
\numberwithin{equation}{section}
\begin{document}

\author{Jakob Scholbach}

\title{Arakelov motivic cohomology II}
%\author{Jakob Scholbach \footnote{Universit{\"a}t M{\"u}nster, Mathematisches Institut, Einsteinstr. 62, D-48149 M{\"u}nster, Germany,
%\href{mailto:jakob.scholbach@uni-muenster.de}{jakob.scholbach@uni-muenster.de}
%}}
%\subclass{14G40, 14F35, 19E08}
%\keywords{Arakelov theory, motives, spectrum, height pairing, arithmetic Chow groups, arithmetic K-groups}

\maketitle

%\tableofcontents

\begin{abstract}
We show that the constructions done in part I generalize their classical counterparts: firstly, the classical Beilinson regulator is induced by the abstract Chern class map from $\BGL$ to the Deligne cohomology spectrum. Secondly, Arakelov motivic cohomology is a generalization of arithmetic $K$-theory and arithmetic Chow groups. For example, this implies a decomposition of higher arithmetic $K$-groups in its Adams eigenspaces. Finally, we give a conceptual explanation of the height pairing: it is the natural pairing of motivic homology and Arakelov motivic cohomology.
\end{abstract}

The purpose of this work is to compare the abstract constructions of the regulator map and the newly minted Arakelov motivic cohomology groups done in part I with their classical, more concrete counterparts. In a nutshell, Arakelov motivic generalizes and simplifies a number of classical constructions pertaining to arithmetic $K$- and Chow groups.

We show that the Chern class $\ch_\Deligne: \BGL \r \oplus_p \HD\twi p$ between the spectra representing $K$-theory and Deligne cohomology constructed in \ref{defi_regulators} induces the Beilinson regulator
$$K_n(X) \r \oplus_p \HD^{2p-n}(X, p)$$
for any smooth scheme $X$ over an arithmetic field (\refth{comparisonRegulator}).

Next, we turn to the relation of Arakelov motivic cohomology and arithmetic $K$- and Chow groups. Arithmetic $K$-groups were defined by Gillet-Soul\'e and generalized to higher $K$-theory by Takeda \cite{GilletSoule:CharacteristicI, GilletSoule:CharacteristicII, Takeda:Higher}. We denote these groups by $\Khat^T_n(X)$. They fit into an exact sequence
$$K_{n+1}(X) \r \Deligne_{n+1}(X) / \im d_\Deligne \r \Khat^T_n(X) \r K_n(X) \r 0,$$
where $\Deligne_*(X)$ is a certain complex of differential forms. The presence of the group $\Deligne_{n+1}(X) / \im d_\Deligne$, as opposed to the Deligne cohomology group $\ker d_\Deligne / \im d_\Deligne = \oplus_p \HD^{2p-n-1}(X, p)$ implies that the groups $\Khat^T_n(X)$ are not homotopy invariant. Therefore they cannot be adressed using $\A$-homotopy theory. Instead, we focus on the subgroup (see p.\ \pageref{eqn_KhatTakeda}),
$$\Khat_n(X) := \ker \left (\ch: \Khat^T_n(X) \r \Deligne_n(X) \right ).$$
and show a canonical isomorphism
\eqn \mylabel{eqn_Khatintro}  \tag{*}
\Hhat^{-n}(X) \cong \Khat_n(X).
\xeqn
for smooth schemes $X$ and $n \geq 0$. All our comparison results concern the groups $\Khat_*(X)$ and, in a similar vein, the subgroup $\CHhat^*(X)$ of Gillet-Soul\'e's group \cite{GilletSoule:Arithmetic} consisting of arithmetic cycles $(Z, g)$ satisfying $\delta_Z = \partial \overline \partial g / (2 \pi i)$, cf.\ p.\ \pageref{eqn_CHhat}. The homotopy-theoretic approach taken in this paper conceptually explains, improves, and generalizes classical constructions such as the arithmetic Riemann-Roch theorem, as far as these smaller groups are concerned. The simplification stems from the fact that it is no longer necessary to construct explicit homotopies between the complexes representing arithmetic $K$-groups, say. For example, the Adams operations on $\Khat_n(X)$ defined by Feliu \cite{Feliu:Adams} were not known to induce a decomposition
$$\Khat_*(X)_\Q \cong \oplus_p \Khat_*(X)^{(p)}_\Q.$$
Using that the isomorphism \refeq{Khatintro} is compatible with Adams operations, this statement follows from the entirely formal analogue for $\Hhat^*$, namely the Arakelov-Chern class isomorphism \refeq{HBhatVsBGLhat}. We conclude a canonical isomorphism
$$\Hhat^{2p,p}(X, p) = \CHhat^p(X)_\Q.$$
Moreover, the pushforward on Arakelov motivic cohomology established in \ref{defilemm_pushforward} is shown to agree with the one on arithmetic Chow groups in two cases, namely for the map $\SpecFp \r \SpecZ$ and for a smooth proper map $X \r S$, $S \subset \SpecOF$ for a number ring $\OF$. The non-formal input in the second statement is the  finiteness of the Chow group $\CH^{\dim X}(X)$ proven by Kato and Saito \cite{KatoSaito:Global}. In a similar vein, we identify the pushforward on $\Khat_0$ with the one on $\Hhat^0$ (\refth{comparisonpushforward}). %In particular, the $\Khat_0$-theoretic pushforward does not depend on the choices of any metrics.
In \refsect{motivicduality}, it is shown that the height pairing
$$\CH^m(X) \x \CHhat^{\dim X-m}(X) \r \CHhat^1(S)$$
coincides, after tensoring with $\Q$, with the \emph{Arakelov intersection pairing} of the motive $M := \M(X)(m-\dim X+1)[2(m-\dim X+1)]$ of any smooth proper scheme $X / S$:
\eqnarr
\Hom_{\SH(S)}(S^0, M) \x \Hom (M, \HBhatS (1)[2]) & \r&  \Hhat^2(S, 1), \\
(\alpha, \beta) & \mapsto & \beta \circ \alpha.
\xeqnarr
Conjecturally, the $L$-values of schemes (or motives) over $\Z$ are given by the determinant of this pairing \cite{Scholbach:SpecialL}.

%Moreover, the higher arithmetic Riemann-Roch theorem for $\Hhat(X)$ vs. $\Hhat^*(X,*)$ (\refth{HARR}) shows that the arithmetic Riemann-Roch theorem of Gillet, Roessler, and Soul\'e is  essentially a consequence of motivic considerations, at least as far as the subgroups $\Khat_0(X) \subset \Khat^T_0(X)$ are concerned. Of course, the approach taken here also provides a conceptual generalization to higher arithmetic $K$-groups.

In the light of these results, stable homotopy theory offers a conceptual clarification of hitherto difficult or cumbersome explicit constructions of chain maps and homotopies representing the expected maps on arithmetic $K$-theory, such as the Adams operations. The bridge between these concrete constructions and the abstract path taken here is provided by a strong unicity theorem. Recall that there is a distinguished triangle
$$\oplus_{p \in \Z} \HD \twi p[-1] \r \BGLhat \r \BGL \stackrel{\ch_\Deligne}\lr \oplus_{p \in \Z} \HD \twi p$$
in the stable homotopy category. Among other things we prove that $\BGLhat$ is unique, up to \emph{unique} isomorphism fitting into the obvious map of distinguished triangles (see \ref{theo_comparisongroups} for the precise statement). The proof of this theorem takes advantage of the motivic machinery, especially the computations of Riou pertaining to endomorphisms of $\BGL$. Its only non-formal input is a mild condition involving the $K$-theory and Deligne cohomology of the base scheme. The unicity trickles down to the unstable homotopy category. It can therefore be paraphrased as: any construction for the groups $\Khat_*$ that is functorially representable by zig-zags of chain maps and compatible with its non-Arakelov counterpart is necesssarily unique. The above-mentioned identification of the Adams operations and the $K$-theory module structure on $\Khat$ are consequences of this principle. In order to show that the arithmetic Riemann-Roch theorem by Gillet, Roessler and Soul\'e \cite{GilletRoesslerSoule}, when restricted to $\Khat_0(X) \subset \Khat^T_0(X)$ (!), is a formal consequence of the motivic framework it remains to show that their arithmetic Chern class \cite[cf. Thm. 7.2.1]{GilletSoule:CharacteristicII},
$$\Khat_0(X)_\Q \cong \oplus_p \Khat_0(X)^{(p)}_\Q,$$
agrees with the Arakelov Chern class established in \refeq{HBhatVsBGLhat}. This will be a consequence of the above unicity result, once the arithmetic Chern class has been extended to higher arithmetic $K$-theory by means of a canonical (i.e., functorial) zig-zag of appropriate chain complexes.
\\

%\end{TheoStar}

\noindent \emph{Acknowledgement: }%I would like to thank J\"urg Kramer for raising the question how Arakelov Chow groups $\CHbar$ are related to Arakelov motivic cohomology. Also,
I would like to thank Andreas Holmstrom for the collaboration leading to part I of this project.

\setcounter{section}{4}

\section{Comparison of the regulator} \mylabel{sect_comparisonRegulator}
After recalling some details of the construction of $\BGL$ in \refsect{prelimBGL}, we construct a Chern class map $\ch: \BGL \r \oplus_p \HD\twi p$ that induces the Beilinson regulator. This is done in \refsect{secondconstruction}. The strategy is to take Burgos' and Wangs representation of the Beilinson regulator as a map of simplicial presheaves and lift it to a map in $\SH(S)$. We finish this section by proving that this Chern class $\ch$  and the one obtained in \refde{regulators},
$$\ch_\Deligne : \BGL \stackrel{\id \wedge 1_\Deligne} \lr \BGL_\Q \wedge \HD \stackrel{\ch \wedge \id}\lr \oplus_{p \in \Z} \HB \twi p \wedge \HD \stackrel {1_{\Beilinson} \wedge \id_\Deligne, \cong} \longleftarrow \oplus_p \HD \twi p,$$
agree. In particular, $\ch_\Deligne$ also induces the Beilinson regulator. This result is certainly not surprising---after all Beilinson's regulator is the Chern character map for Deligne cohomology.

Throughout, we will use the notation of part I. In particular, $\Ho_\bullet(S)$ and $\SH(S)$ are the unstable and the stable homotopy category of smooth schemes over some Noetherian base scheme $S$ (Sections \ref{sect_SH}, \ref{sect_motives}).

\subsection{Reminders on the object $\BGL$ representing $K$-theory} \mylabel{sect_prelimBGL}

In order to prove our comparison results, we need some more details concerning the object $\BGL$ representing algebraic $K$-theory. This is due to Riou \cite{Riou:Thesis}.

Let $\Gr_{d,r}$ be the Grassmannian whose $T$-points, for any $T \in \Sm /S$, are given by locally free subsheaves of $\cO^{d+r}_T$ of rank $d$. As usual, we regard this (smooth projective) scheme as a presheaf on $\Sm / S$. For $d \leq d'$, $r \leq r'$, the transition map
\eqn \mylabel{eqn_transitionGrass}
\Gr_{d,r} \r \Gr_{d',r'}
\xeqn
is given on the level of $T$-points by mapping $M \subset \cO_T^{d+r}$ to $\cO_T^{d'-d} \oplus M \oplus 0^{r'-r} \subset \mathcal O^{d'+r'}$. Put $\Gr := \varinjlim_{\N^2} \Gr_{*, *}$, where the colimit is taken in $\PSh(\Sm / S)$. It is pointed by the image of $\Gr_{0,0}$. Write $\Z \x \Gr$ for the product of the constant sheaf $\Z$ (pointed by zero) and this presheaf, and also for its image in $\HoC_\bullet(S)$. For a regular base scheme $S$, there is a functorial (with respect to pullback) isomorphism
\eqn \mylabel{eqn_Krepunstable}
\Hom_{\Ho_\bullet(S)}(S^n \wedge X_+, \Z \x \Gr) \cong K_n(X),
\xeqn
for any $X \in \Sm/S$ \cite[Prop. 3.7, 3.9, page 138]{MorelVoevodsky:A1}.

\defi \cite[I.124, IV.3]{Riou:Thesis} \mylabel{defi_SHnaive}
%\begin{enumerate}[(i)]
%\item
The category $\SHnaive(S)$ is the category of $\Omega$-spectra (with respect to $-\wedge \P$) in $\HoC_\bullet(S)$: its objects are sequences $E_n \in \HoC_\bullet(S)$, $n \in \N$ with bonding maps $\P \wedge E_n \r E_{n+1}$ inducing isomorphisms $E_n \r \IHom_\bullet(\P, E_{n+1})$.\footnote{We will not write $\mathrm L$ or $\R$ for derived functors. For example, $f^*$ stands for what is often denoted $\mathrm L f^*$ and similarly with right derived functors such as $\R \IHom$, $\R \Omega$ etc.}  %(total derived inner $\Hom$ of pointed sheaves).
Its morphisms are sequences of maps $f_n : E_n \r F_n$  (in $\HoC_\bullet(S)$) making the diagrams involving the bonding maps commute.
%\end{enumerate}
\xdefi

\rema \mylabel{rema_SHnaiveforgetful}
Recall the \emph{projective Nisnevich-$\A$-model structure} on $\P$-spectra: a map $f: X \r Y$ is a weak equivalence (fibration), if all its levels $f_n : X_n \r Y_n$ is a weak equivalence (fibration, respectively) in the Nisnevich-$\A$-model structure on $\Delta^\op(\PSh_\bullet(\Sm/S))$ (whose homotopy category is $\HoC_\bullet(S)$. The homotopy category of spectra with respect to the projective model structure is denoted $\SH_p(S)$. The composition of the inclusion of the full subcategory of $\Omega$-spectra and the natural localization functor,
$$\{X \in \SH_p, X \text{ is an } \Omega-\text{spectrum}\} \subset \SH_p(S) \r \SH(S),$$
is an equivalence. %This equivalence is a consequence of \cite[2.8, 4.6]{Jardine:Motivic}.
This yields a natural ``forgetful'' functor $\SH(S) \r \SHnaive(S)$.
\xrema

\defitheo (Riou, \cite[IV.46, IV.72]{Riou:Thesis}) \mylabel{defitheo_BGLnaive} %\mylabel{theo_BGL}
The spectrum $\BGLnaive \in \SHnaive(S)$ consists of $\BGLnaive_n := \Z \x \Gr$ (for each $n \geq 0$) with bonding maps
\eqn \mylabel{eqn_bondingZGr}
\P \wedge (\Z \x \Gr) \stackrel{u_1^* \wedge \id}\lr (\Z \x \Gr) \wedge (\Z \x \Gr) \stackrel \mu \lr \Z \x \Gr,
\xeqn
where $u_1^*$ is the map corresponding to $u_1 = [\cO(1)] - [\cO(0)] \in K_0 (\P) \stackrel{\refeq{Krepunstable}}= \Hom_{\Ho}(\P, \Z \x \Gr)$ and $\mu$ is the multiplication map, that is to say, the unique map  \cite[III.31]{Riou:Thesis} inducing the natural (i.e., tensor) product on $K_0(-)$.

For $S = \SpecZ$, there is a lift $\BGL_\Z \in \SH(\SpecZ)$ of $\BGLnaive \in \SHnaive(\Z)$ that is unique up to \emph{unique} isomorphism. For any scheme $f: S \r \SpecZ$, put $\BGL_S := f^* \BGL_\Z$. The unstable representability theorem \refeq{Krepunstable} extends to an isomorphism
\eqn \mylabel{eqn_BGL_vs_SH_nochmal}
\Hom_{\SH(S)}(S^n \wedge \Sigma_{\P}^\infty X_+, \BGL_S) = K_n(X)
\xeqn
for any regular scheme $S$ and any smooth scheme $X / S$. In $\SH(S)_\Q$, i.e., with rational coefficients, $\BGL_S \t \Q$ decomposes as
\eqn \mylabel{eqn_BGLdecomposition}
\BGL_S \t \Q = \oplus_{p \in \Z} \HBS(p)[2p]
\xeqn
such that the pieces $\HBS(p)[2p]$ represent the graded pieces of the $\gamma$-filtration on $K$-theory:
$$ %\eqn
%\mylabel{eqn_motivic_vs_K_nochmal}
\Hom_{\SH(S)}(S^n \wedge \Sigma_{\P}^\infty X_+, \HBS(p)[2p]) \cong \gr_\gamma^p K_n(X)_\Q.
$$ %\xeqn
\xdefitheo
%\pf
%This is due to Riou. Previously, Morel and Voevodsky showed that $K$-theory is representable in $\HoC_\bullet(S)$ \cite[Th. 3.13, p. 140]{MorelVoevodsky:A1}, and Riou was able to lift the Adams operations on $K$-theory to the representing object.
%\xpf

\lemm \mylabel{lemm_Grassmannian}
For any $d$, $r$, the motive $\M(\Gr_{d,r})$ (cf. \refsect{motives}) is given by
\eqn
\mylabel{eqn_Grassmannian}
\M(\Gr_{d,r}) = \oplus_\sigma \M(S) \left (\sum (\sigma_i - i) \right) \left [2\sum (\sigma_i - i) \right].
\xeqn
The sum runs over all Schubert symbols, i.e., sequences of integers satisfying $1 \leq \sigma_1 < \dots < \sigma_d \leq d+r$.
For $d \leq d'$, $r \leq r'$, the transition maps \refeq{transitionGrass} $\M(\Gr_{d,r}) \r \M(\Gr_{d',r'})$ exhibits the former motive as a direct summand of the latter.
\xlemm
\pf
Formula \refeq{Grassmannian} is well-known \cite[2.4]{Semenov:Motives}. The second statement follows from the same technique, namely the localization triangles for motives with compact support applied to the cell decomposition of the Grassmannian:
%(Hatcher, Vector bundles and K-theory, p. 32
for any field $k$, a $d$-space $V^{(d)}$ in $k^{d+r}$ is uniquely described by a $(d, d+r)$-matrix $A$ in echelon form %(the matrix corresponds the image of the linear map $k^d \r k^{d+r}$ it describes)
such that $A_{\sigma_i,j} = \delta_{i,j}$ and $A_{i, j} = 0$ for $i > \sigma_j$ for some Schubert symbol $\sigma$. The constructible subscheme of $\Gr_{d,r}$ whose $k$-points are given by matrices with fixed $\sigma$ is an affine space $\A[(\sigma)]_S$. % of dimension $\sum (\sigma_i - i)$.
The transition map $V^{(d)} \mapsto k^{d'-d} \oplus V^{(d)} \oplus 0^{r'-r}$ corresponds to
$$A \mapsto \left [ \begin{array}{ccc}
\Id_{d'-d} & 0 & 0 \\
       0 & A & 0 \\
0 & 0 & 0^{r'-r}
\end{array}
\right ],
$$
that is,
$$\sigma \mapsto \left(1, 2, \dots, d'-d, \sigma_1+(d'-d), \dots \sigma_d+(d'-d)\right ) =: \sigma'.$$
In other words, the restriction of the transition maps \refeq{transitionGrass} to the cells is the identity map $\A[(\sigma)]_S \r \A[(\sigma')]_S$, which shows the second statement.
\xpf

\subsection{Second construction of the regulator} \mylabel{sect_secondconstruction}
In this subsection and the next one, $S$ is an arithmetic field and $X$ is a smooth scheme over $S$.

Let
%\eqn \mylabel{eqn_DoldKan}
$\cK: \Com^{\geq 0}(\Ab) \r \Delta^\op \Ab$ %: \mathcal N
%\xeqn
be the Dold-Kan equivalence on chain complexes concentrated in degrees $\geq 0$ (with $\deg d = -1$ and shift given by $C[-1]_a = C_{a-1}$). Recall from Definitions \ref{defi_Burgos}, \ref{defi_Delignesimplicial} the abelian presheaf complex $\Deligne$ and $\Deligne_s := \cK(\tau_{\geq 0} \Deligne)$. We have $\H_n(\Deligne(X)) = \pi_n(\Deligne_s(X)) = \oplus_p \HD^{2p-n}(X, p)$. We set $\Deligne_s[-1] := \cK ((\tau_{\geq 0} \Deligne)[-1])$. Recall that for any chain complex of abelian groups $C$, %with degree of differential is -1
there is a natural map $S^1 \wedge \cK(C) = \cone (\cK(C) \r \text{point}) \r \cK(\cone(C \r 0)) = \cK(C[-1])$, hence a map
%\eqn \mylabel{eqn_ROmegaDoldKan}
$\cK(C) \r \Omega_s \cK(C[-1])$.
%\xeqn
(Here and elsewhere, $\Omega_s$ is the simplicial loop space, its $\P$-analogue is denoted $\Omega_{\P}$.) This map is a weak equivalence of simplicial abelian groups.

For any pointed simplicial presheaf $F \in \HoC_\bullet(S)$, let $\varphi(F)$ be the pointed presheaf
\eqn
\mylabel{eqn_functorphi}
\varphi(F): \Sm / S \ni X \mapsto \Hom_{\HoC_\bullet(S)}(X_+, F).
\xeqn
According to \refeq{Krepunstable} and \refle{DeligneHoC}, respectively,
\eqnarra
\varphi(\Z \x \Gr) & =&  K_0: X \mapsto K_0(X), \mylabel{eqn_varphiK0}\\
\varphi(\Omega_s^n \Deligne_s) & = & \HD^{-n}: X \mapsto \oplus_p \HD^{2p-n}(X, p), \ n \geq 0. \nonumber
\xeqnarra

\mylabel{hatP} Let $\Phat(X)$ be the (essentially small) Waldhausen category consisting of hermitian bundles $\overline E = (E, h)$ on $X$, i.e., a vector bundle $E / X$ with a metric $h$ on $E(\CC) / X(\CC)$ that is invariant under $\Fr_\infty^*$ and smooth at infinity \cite[Definition 2.5]{BurgosWang}. Morphisms are given by usual morphisms of bundles, ignoring the metric, so that $\Phat(X)$ is equivalent to the usual category of vector bundles. Let
\eqn
\mylabel{eqn_Sconstruction}
S_*: \Sm / S \ni X \mapsto \text{Sing} |S_* \Phat(X)|
\xeqn
be the presheaf (pointed by the zero bundle) whose sections are given by the simplicial set of singular chains in the topological realization of the Waldhausen $S$-construction of $\Phat(X)$. %see e.g.\ \cite[1.4]{Carlsson:Deloopings} for the $S$-construction
Its homotopy presheaves are
\eqn \mylabel{eqn_SconstructionVsK}
\Hom_{\HoCsect(S)} (S^n \wedge X_+, S_*) = \pi_n S_*(X) = \pi_{n-1} \Omega_s S_*(X) \cong K_{n-1}(X), \ \ n \geq 1.
\xeqn
Here, $\HoCsect$ denotes the homotopy category of $\Delta^\op \PSh_\bullet(\Sm/S)$ (simplicial pointed presheaves), endowed with the section-wise model structure. $K$-theory (of regular schemes) is homotopy invariant and satisfies Nisnevich descent \cite[Thm. 10.8]{ThomasonTrobaugh}. Therefore, as is well-known, the left hand term agrees with $\Hom_{\Ho_\bullet(S)}(S^n \wedge X_+, S_*)$. That is, there is an isomorphism of pointed presheaves
\eqn \mylabel{eqn_OmegaSK0}
\varphi (\Omega_s S_*) \cong K_0.
\xeqn
According to \cite[III.61]{Riou:Thesis}, there is a unique isomorphism in $\HoC_\bullet(S)$
\eqn \mylabel{eqn_isoGrWaldhausen}
\tau: \Z \x \Gr \r \Omega_s S_* .
\xeqn
making the obvious triangle involving \refeq{OmegaSK0} and \refeq{varphiK0} commute.
%This relies on the interpretation of $K$-theory in terms of morphisms in $\HoC_\bullet(S)$ \cite[3.9., p. 139]{MorelVoevodsky:A1}.

The proof of our comparison of the regulator uses the following result due to Burgos and Wang \cite[Prop. 3.11, Theorem 5.2., Prop. 6.13]{BurgosWang}:
\satz \mylabel{prop_BurgosWang}
There is a map in $\Delta^\op(\PSh_\bullet(\Sm/S))$
$$\ch_S: S_* \r \Deligne_s[-1]$$
such that the induced map
\eqnarr
\pi_{n} \ch_S: K_{n-1}(X) \r \oplus_{p \in \Z} \HD^{2p-(n-1)}(X, p) %& = & \pi_{n}(S_*(X)) \\
\xeqnarr
agrees with the Beilinson regulator for all $n \geq 1$.
\xsatz

By \refeq{isoGrWaldhausen}, %and \refeq{ROmegaDoldKan}
we get a map in $\HoC_\bullet(S)$:
\eqn \mylabel{eqn_ChernHoC}
\ch: \Z \x \Gr \stackrel {\tau, \cong} \lr \Omega_s S_* \stackrel{\Omega_s \ch_S}\lr \Omega_s (\Deligne_s[-1]) \stackrel \cong \lr \Deligne_s.
\xeqn
The induced map
\eqn \mylabel{eqn_regulatorDownToEarth}
K_n(X) \stackrel{\text{\refeq{BGL_vs_SH_nochmal}}}\cong \Hom_{\HoC_\bullet}(S^n \wedge X_+, \Z \x \Gr) \r \Hom_{\HoC_\bullet}(S^n \wedge X_+, \Deligne_s)
\stackrel{\text{\refeq{DeligneHoC}}}\cong \oplus_p \HD^{2p-n}(X, p)
\xeqn
agrees with the Beilinson regulator. %, at least up to the identification of $K$-theory as represented by $\Z \x \Gr$ vs.\ as computed by the $S$-construction, cf.\ \refeq{BGL_vs_SH}, \refeq{SconstructionVsK}. For the groups $K_0(X)$, this identification is the usual one, mapping the class $[E]$ of some vector bundle to the corresponding $1$-simplex in $S(X)$ \cite[III.61]{Riou:Thesis}. Possible applications of Arakelov motivic cohomology to special $L$-values of motives are concerned with the covolume of maps such as \refeq{regulatorDownToEarth}, which is independent of the concrete complex whose homotopy (or homology) groups compute $K$-theory. Therefore, we are going to ignore this issue in this paper.
In order to lift the map $\ch$ to a map in $\SH(S)$, we first check the compatibility with the $\P$-spectrum structures to get a map in $\SHnaive(S)$. This means that the diagram involving the bonding maps only has to commute up to ($\A$-)homotopy. Then, we apply an argument of Riou to show that this map actually lifts uniquely to one in $\SH(S)$. %Let us point out that the theorem below is a close adaptation of Riou's work and could also be derived as a corollary of the analogous statements concerning the regulator map taking values in motivic cohomoloy.

Recall the Deligne cohomology ($\P$-)spectrum $\cDz$ from \refdele{Delignespectrum}. Its $p$-th level is given by $\Deligne_s(p)$, for any $p \geq 0$.

\theo \mylabel{theo_mapSH}
\begin{enumerate}[(i)]
\item \mylabel{item_pfmapSHi}
In $\SHnaive(S)$, there is a unique map
$$\ch^\mathrm{naive}: \BGLnaive_S \r \oplus_{p \in \Z}\cDz(p)[2p] =: \oplus_p \HD \twi p$$
that is given by $\ch: \Z \x \Gr \stackrel{\text{\refeq{ChernHoC}}}\lr \Deligne_s$ in each level.
\item \mylabel{item_pfmapSHii}
In $\SH(S)$, there is a unique map
$$\ch: \BGL_S \r \oplus_{p \in \Z}\cDz(p)[2p]$$
that maps to $\ch^\mathrm{naive}$  under the forgetful functor $\SH(S) \r \SHnaive(S)$ (\refre{SHnaiveforgetful}).
\item \mylabel{item_pfmapSHiii}
There is a unique map
$$\rho: \HBS \r \cDz$$
in $\SH(S)_\Q$ such that $\ch \t \Q = \oplus_{p \in \Z} \rho(p)[2p]: \BGL_\Q \r \oplus \cDz(p)[2p]$, under the identification \refeq{BGLdecomposition}.
\end{enumerate}
\xtheo
\pf
By \refle{Grassmannian}, the transition maps \refeq{transitionGrass} defining the infinite Grassmannian induce split monomorphisms $\M(\Gr_{d,r}) \r \M(\Gr_{d',r'})$ of motives and therefore (e.g.\ using \refth{Delignespectrum}) split surjections (for any $n \geq 0$, $d \leq d'$, $r \leq r'$)
\eqn
\mylabel{eqn_surjectivityGrassmannian}
\begin{array}{ccc}
\Hom_{\HoC(S)} (\Gr_{d', r'}, \Omega_s^n \Deligne_s) & \r & \Hom_{\HoC(S)} (\Gr_{d, r}, \Omega_s^n \Deligne_s) \\
\Vert & & \Vert \\
\HD^{-n}(\Gr_{d',r'}) & & \HD^{-n}(\Gr_{d,r}).
\end{array}
\xeqn
%(This causes a certain $\R^1 \varprojlim$ to vanish.)
%Indeed $\HD^*(\Gr_{d,r}, *) = \Hw^*(\underline \RG_{\MHS}(\Gr_{d,r})(*))$. Here $\Hw^*(-)$ is the weak Hodge cohomology of a complex of mixed Hodge structures \todo{REF BEILINSON} and $\underline \RG_{\MHS}(-)$ is the complex of Hodge structures whose cohomology objects are the cohomology of the variety in question, with its Hodge structure as established by Deligne. By \refle{Grassmannian}, \refeq{transitionGrass} induce split monomorphisms $\M(\Gr_{d,r}) \r \M(\Gr_{d',r'})$ of motives, thus split epimorphisms after applying $\underline \RG_{\MHS}$, hence for the Deligne cohomology.
A similar surjectivity statement holds for the map of Deligne cohomology groups induced by transition maps defining the product $\Gr \x \Gr$, i.e.,
$$\Gr_{d_1, r_1} \x \Gr_{d_2, r_2} \r \Gr_{d'_1, r'_1} \x \Gr_{d'_2, r'_2}.$$

\refit{pfmapSHi}: the unicity of $\ch^\mathrm{naive}$ is obvious. Its existence amounts to the commutativity of the following diagram in $\HoC_\bullet(S)$:
\eqn
\mylabel{eqn_SHnaive}
\xymatrix{
\P \wedge \Z \x \Gr \ar[d]^{\id \wedge \ch} \ar[r]^<<<<{u_1^* \wedge \id} &
(\Z \x \Gr) \wedge (\Z \x \Gr) \ar[r]^>>>>>{\mu} \ar[d]^{\ch \wedge \ch} &
\Z \x \Gr \ar[d]^{\ch} \\
\P \wedge \Deligne_s  \ar[r]^{c^* \wedge \id}
 &
\Deligne_s \wedge \Deligne_s \ar[r]^{\mu}
 &
\Deligne_s.
}
\xeqn
The top and bottom lines are the bonding maps of $\BGLnaive$ (cf. \refeq{bondingZGr}) and $\oplus_p \cDz \twi p$ (cf. \refdele{Delignespectrum}), respectively. The map $c^*$ corresponds to the first Chern class $c_1(\cO_{\P}(1)) \in \HD^2(\P_S, 1)$. To see the commutativity of the right half, we use that the functor $\varphi$ \refeq{functorphi} induces an isomorphism
$$\Hom_{\HoC_\bullet(S)}((\Z \x \Gr)^{\wedge 2}, \Deligne_s) = \Hom_{\PSh_\bullet(\Sm/S)}(K_0(-) \wedge K_0, \HD^0).$$
This identification is shown exactly as \cite[III.31]{Riou:Thesis}, which does it for $\Z \x \Gr$ instead of $\Deligne_s$. The point is a surjectivity argument in comparing cohomology groups of products of different Grassmannianns, which is applicable to Deligne cohomology by the remark above. By construction of the multiplication map on $\Z \x \Gr$, applying $\varphi$ to the right half of \refeq{SHnaive} yields the diagram
$$\xymatrix{
K_0 \wedge K_0 \ar[r]^{\mu_{K_0}} \ar[d]^{\ch \wedge \ch} &
K_0 \ar[d]^\ch \\
\HD^0 \wedge \HD^0 \ar[r]^{\mu_{\Deligne}} &
\HD^0.
}$$
Here $\mu_{K_0}$ is the usual (tensor) product on $K_0$ and $\mu_\Deligne$ is the classical product on Deligne cohomology \cite{EsnaultViehweg:Deligne}. The Beilinson regulator is multiplicative \cite[Cor., p. 28]{Schneider:Introduction}, so this diagram commutes.

%It is sufficient to show that $\mu \circ (\ch \x \ch) = \ch \circ \mu$ (as opposed to $\wedge$ in place of $\x$).
For the commutativity of the left half, let $i_{m,n}: \P[m] \r \P[n]$ be the inclusion $[x_0: \hdots : x_m] \mapsto [x_0: \hdots : x_m : 0 : \hdots : 0]$, for $m \leq n$, and $i_{m, \infty} := \colim_n i_{m,n} : \P[m] \r \P[\infty] := \colim_n \P[n]$. The map $u_1^*$ factors as
$$\P \stackrel{i_{1, \infty}}\lr \P[\infty] \stackrel{u_\infty^*} \lr \Z \x \Gr$$
where $u_\infty^* \in \Hom_{\Ho_\bullet(S)}(\P[\infty], \Z \x \Gr)$ is induced by the compatible system $u_n = [\cO_{\P[n]}(1)] - [\cO_{\P[n]}] \in K_0(\P[n])$---simply because $i_{1,n}^* \cO_{\P[n]}(1) = \cO_{\P}(1)$. Similarly, $c^* = c_1(\cO(1))$ is given by
$$c^*: \P \stackrel{i_{1, \infty}} \lr \P[\infty] \stackrel{u_\infty^*}\lr \Z \x \Gr \stackrel{\ch}\lr \Deligne_s,$$
%factors over the map $\P[\infty] \r \Deligne_s$ (in $\HoC_\bullet(S)$) given by the compatible system of $c_1(\cO_{\P[n]}(1))$ which in turn  factors as
%$$\P[\infty] \stackrel{u_\infty^*}\lr \Z \x \Gr \stackrel{\ch}\lr \Deligne_s,$$
because $\ch(\cO(1)) - \ch(\cO) = \exp(c_1(\cO(1)) - 1$ which on $\P$ equals $c_1(\cO(1)) \in \HD^{2}(\P[\infty], 1)$. Then the commutativity of the diagram in question is obvious.
%\scriptsize
%Therefore, the commutativity of the upper half is a consequence of the obviously commutative diagram (the horizontal maps are identity wedge $\ch$):
%
%$$\xymatrix{
%\P[\infty] \wedge (\Z \x \Gr) \ar[r] \ar[d]^{u_\infty^* \wedge \id} &
%\P[\infty] \wedge \cD_n \ar[d]^{u_\infty^* \wedge \id} \\
%(\Z \x \Gr) \wedge (\Z \x \Gr) \ar[r] \ar[d]^{\ch \wedge \id} &
%(\Z \x \Gr) \wedge \cD_n \ar[d]^{\ch \wedge \id} \\
%\cD_1 \wedge (\Z \x \Gr) \ar[r] &
%\cD_1 \wedge \cD_n
%}
%$$
%\normalsize

\refit{pfmapSHii}: For each $n \geq 0$ and $m =0, -1$, put $V^m_n := \Hom_{\PSh(\Sm/S, \Ab)}(K_0, \HD^m)$. These groups form a projective system with transition maps
$$V^m_{n+1} \ni (f_n : K_0 \r \HD^m) \mapsto (\Omega_{\P} f_n : \Omega_{\P} K_0 \r \Omega_{\P} \HD^m) \in V^m_n,$$
where $\Omega_{\P} (F)$ is the presheaf $\Sm / S \ni U \mapsto \ker (F(\P_U) \stackrel{\infty^*}\r F(U))$. Indeed, the projective bundle formula (for $\P$) implies an isomorphism of presheaves $\Omega_{\P} K_0 \cong K_0$ and likewise with $\HD^m$.

The composition of functors
$$\SH \r \SH^{\text{naive}} \stackrel n \lr \Ho_\bullet \stackrel \varphi \lr \PSh (\Sm / S)$$
actually takes values in $\PSh(\Sm/S, \Ab)$. Here, $n$ indicates taking the $n$-th level of a spectrum. By construction, $\BGL$ gets mapped to $K_0$ and $\HD$ gets mapped to the presheaf $\HD^0 = \oplus_{p} \HD^{2p}(-, p)$ for each $n \geq 0$. This gives rise to the following map (cf. \cite[IV.11]{Riou:Thesis})
$$\Hom_{\SH}(\BGL, \oplus_{p}\cDz \twi p) \r \Hom_{\SHnaive(S)}(\BGLnaive, \oplus_{p}\cDz \twi p) \cong \varprojlim_n V^0_n.$$
This map is part of the following Milnor-type short exact sequence \cite[IV.48, III.26, see also IV.8]{Riou:Thesis} (it is applicable because of the surjectivity of \refeq{surjectivityGrassmannian} for $n=1$ and $n=2$)
\eqn
\mylabel{eqn_Milnorexact}
0 \r \R^1 \varprojlim V^{-1}_n \r \Hom_{\SH}(\BGL, \oplus_{p}\cDz \twi p) \r \varprojlim_n V^0_n \r 0.
\xeqn
The map $\ch^{\mathrm{naive}}$ thus corresponds to a unique element in the right-most term of \refeq{Milnorexact}.
The natural map
\small
\eqnarr
V^{-1}_n = \Hom_{\PSh(\Ab)}(K_0, \HD^{-1}) & \r & \varprojlim_e \oplus_p \HD^{2p-1}(\P[e]_S,p) \cong \oplus_{p \in \Z} \oplus_{j=0}^p \HD^{2p-2j-1}(S, p-j)\\
f & \mapsto & \left ( f(\cO_{\P[e]}(1)) \right)_e,
\xeqnarr
\normalsize
is an isomorphism. Indeed, the proof of the analogous statement for motivic cohomology instead of Deligne cohomology \cite[V.18]{Riou:Thesis} (essentially a splitting argument) only uses the calculation of motivic cohomology of $\P[e]$. Thus it goes through by the projective bundle formula for Deligne cohomology. %\refeq{projectiveBundleDeligne}.

Via this identification, the transition maps $\Omega_{\P}: V^{-1}_{n+1} \r V^{-1}_n$ are the direct sum over $p \in \Z$ of the maps
$$\oplus_{j = 0}^p \HD^{2p-2j-1}(S, p-j) \r \oplus_{j = 0}^{p-1} \HD^{2(p-1)-2j-1}(S, (p-1)-j)$$
which are the multiplication by $j$ on the $j$-th summand at the left. Again, this is analogous to \cite[V.24]{Riou:Thesis}. In particular $\Omega_{\P}$ is onto, since Deligne cohomology groups are divisible. Therefore $\R^1 \varprojlim V^{-1}_n = 0,$ so \refit{pfmapSHii} is shown.

\refit{pfmapSHiii}: as in \cite[V.36]{Riou:Thesis}, one sees that $\ch \t \Q$ factors over the projections $\BGL_\Q \r \HB$ and $\oplus_{p \in \Z}\cDz(p)[2p] \r \cDz$.
\xpf

\subsection{Comparison}

\theo \mylabel{theo_comparisonRegulator}
The regulator maps $\ch$, $\rho$ constructed in \refth{mapSH} and the regulator maps $\ch_\Deligne$, $\rho_\Deligne$ obtained in \refde{regulators} agree:
\eqnarr
\ch_\Deligne & = & \ch \in \Hom_{\SH(S)}(\BGL, \oplus_p \HD\twi p), \\
\rho_\Deligne & = & \rho \in \Hom_{\SH(S)_\Q}(\HB, \HD).
\xeqnarr
In particular, $\ch_\Deligne$ also induces the Beilinson regulator $K_n (X) \r \oplus_p \HD^{2p-n}(X, p)$ for any $X \in \Sm / S$, $n \geq 0$.
\xtheo
\pf
The map $\ch$ is a map of ring spectra (i.e., monoid objects in $\SH(S)$): the multiplicativity, i.e.,
$\ch \circ \mu_{\BGL} = \mu_\Deligne \circ (\ch \wedge \ch)$ follows from the right half of the diagram \refeq{SHnaive}. The unitality boils down to $\ch(\mathcal O) = 1 \in \HD^0(S, 0)$. We define an $\BGL$-module structure on $\cD := \oplus_{p \in \Z}\cDz \twi p$ in the usual manner:
$$\BGL \wedge \cD \stackrel{\ch \wedge \id}\lr \cD \wedge \cD \stackrel{\mu}\lr \cD.$$
It is indeed an $\BGL$-module, as one sees using that $\ch$ is a ring morphism.
%\scriptsize
%$$\xymatrix{
%\BGL \wedge \BGL \wedge \cD \ar[dd]^{\mu_{\BGL} \wedge \id_\cD} \ar[r]^{\id \wedge \ch \wedge \id}
%\ar[dr]^{\ch \wedge \ch \wedge \id_\cD}&
%\BGL \wedge \cD \wedge \cD \ar[r]^{\id \wedge \mu_\cD} &
%\BGL \wedge \cD \ar[d]^{\ch \wedge \id_\cD} \\
%&
%\cD \wedge \cD \wedge \cD \ar[d]^{\mu_\cD \wedge \id} &
%\cD \wedge \cD \ar[d]^{\mu_\cD} \\
%\BGL \wedge \cD \ar[r]^{\ch \wedge \id} &
%\cD \wedge \cD \ar[r]^{\mu_\cD} &
%\cD
%}
%$$
%The left part of that diagram uses the multiplicativity of $\ch$.
%The diagram that checks that the unit section of $\BGL$ acts as the identity on $\cD$ is proven similarly.
%\normalsize
%
By the unicity of the $\BGL$-algebra structure on $\cD$ (\refth{Delignespectrum}), this algebra structure agrees with the one established in \lcs This implies $\ch = \ch_\Deligne$. The proof for $\rho = \rho_\Deligne$ is the similar, replacing $\BGL$ with $\HB$ throughout.
%*by \refth{modulestructureDeligne}, .
%In other words, we have shown that the ``universal'' algebra structure on $\cD$ yields the Beilinson regulator.
\xpf

\section{Comparison with arithmetic $K$-theory and arithmetic Chow groups} \mylabel{sect_comparisonKtheory}
In this section, we show that the groups represented by $\BGLhat$ coincide with a certain subgroup of arithmetic $K$-theory as defined by Gillet-Soul\'e and Takeda for smooth schemes over appropriate bases (\refth{comparisongroups}). This isomorphism is compatible with the Adams operations on both sides and with the module structure over $K$-theory (\ref{coro_ChhatVsArakelov}, \ref{theo_comparisonaction}). We also establish the compatibility of the comparison isomorphism with the pushforward in two cases (\ref{theo_comparisonpushforward}).

We consider the following situation: %The geometric set-up in this section is as follows:
%\eqn \mylabel{eqn_geometricsetup}
$X \r S \r B$, where
%\xeqn
$B$ is a fixed arithmetic ring (\refde{arithmring}), $S$ is a regular scheme (of finite type) over $B$ (including the important case $S=B$), and $X \in \Sm/S$. Let $\eta: B_\eta := B \x_\Z \Q \r B$ be the ``generic fiber''. For any datum $?$ related to Deligne cohomology, we write $? := \eta_* ?$ for simplicity of notation. That is, $\Deligne_s(X) := \eta_* \Deligne_s(X) = \Deligne_s(X \x_B B_\eta)$, $\HD := \eta_* \HD \in \SH(S)$ etc.

For a proper arithmetic variety $X$ (i.e., $X$ is regular and flat over an arithmetic ring $B$), Gillet and Soul\'e have defined the arithmetic $K$-group as the free abelian group generated by pairs $(\overline E, \alpha)$, where $\overline E/X$ is a hermitian vector bundle and $\alpha \in \Deligne_0(X) / \im d_{\Deligne}$, modulo the relation
$$(\overline E', \alpha') + (\overline E'', \alpha'') = (\overline E, \alpha'+\alpha''+\widetilde{\ch} (\overline {\mathcal E}))$$
for any extension
$$\overline {\mathcal E}: 0 \r \overline E' \r \overline E \r \overline E'' \r 0$$
of hermitian bundles. Here $\widetilde{\ch}(\overline {\mathcal E})$ is a secondary Chern class of the extension (see \cite[Section 6]{GilletSoule:CharacteristicII} for details). We denote this group by $\Khat^T_0(X)$. The superscript $T$ stands for Takeda, who generalized this to higher $n$ \cite[p. 621]{Takeda:Higher}\footnote{\mylabel{foot_factorsChern} Gillet and Soul\'e use a slightly different normalization  of the Chern class which differs from the one used by Burgos-Wang, Takeda (and this paper) by a factor of $2 (2 \pi i)^n$ for appropriate $n$. See \lcs for details.}. %(see also and Feliu \cite[5.2.2.]{Feliu}).
These higher arithmetic $K$-groups $\Khat^T_n(X)$ fit into a commutative diagram with exact rows and columns, where $\Khat_n (X) := \ker \ch^T$ and $B^\Deligne_n (X) := \im d_\Deligne: \Deligne_{n+1}(X) \r \Deligne_n(X)$:
\small
\eqn \mylabel{eqn_KhatTakeda}
\xymatrix{
%&
%0 \ar[d] &
%0 \ar[d]
%\\
K_{n+1} \ar@{=}[d] \ar[r] &
%\H_{n+1}(\Deligne(X)) =
\oplus_p \HD^{2p-n-1}(p) \ar@{>->}[d] \ar[r]  & \Khat_n \ar@{>->}[d] \ar[r] &  K_n \ar@{=}[d] \ar[r]^(.4){\ch} & \oplus_p \HD^{2p-n}(p)  \\
K_{n+1}  \ar[r] &   \Deligne_{n+1}(X) / \im ({d_\Deligne}) \ar@{->>}[d]^{d_\Deligne} \ar[r] &  \Khat_n^T  \ar@{->>}[d]^{\ch^T} \ar[r] &  K_n  \ar[r] &  0 \\
& B^\Deligne_n %\ar[d]
\ar@{=}[r] & B^\Deligne_n (X) %\ar[d]
%\\
%& 0 & 0
}
\xeqn
\normalsize
The full arithmetic $K$-groups $\Khat^T_*$ are not accessible to homotopy theory since they fail to be $\A$-invariant. Moreover, due to the presence of $\Deligne_{n+1} / \im d_\Deligne$ the groups are usually very large. Therefore, we focus on the subgroups $\Khat_* \subset \Khat^T_*$ and refer to them as arithmetic $K$-theory.

By \refth{comparisonRegulator}, the top exact sequence looks exactly like the one in \refth{longexactArakelov}. In order to show that $\Khat_n(X)$ and $\Hhat^{-n}(X)$ are isomorphic, we use that there is a natural isomorphism (functorial with respect to pullback),
\eqn
\mylabel{eqn_comparisonTakeda}
\Khat_n(X) \cong \pi_{n+1} (\hofib_{\Delta^\op \Sets_\bullet} S_*(X) \stackrel{\ch_S}\lr \Deligne_s[-1](X)), \ \ n \geq 0,
\xeqn
of the arithmetic $K$-group with the homotopy fiber in pointed simplicial sets (endowed with its standard model structure) \cite[Cor. 4.9]{Takeda:Higher}. We write
$$\Shat := \hofib_{\Delta^\op \PSh_\bullet(\Sm/S)} (S_* \r \Deligne_s[-1]),$$
for the homotopy fiber with respect to the section-wise model structure, so that $\pi_{n+1} \left (\Shat(X) \right) = \Khat_n(X)$.

Recall from \refsect{definition} the object $\BGLhat$. Its key property is the existence of a distinguished triangle (in $\SH(S)$)
\eqn \mylabel{eqn_keyproperty}
\oplus_p \HD\twi p[-1] \r \BGLhat \r \BGL \stackrel \ch \r \oplus_p \HD \twi p.
\xeqn
The cohomology groups represented by this object are denoted by $\Hhat^*(-)$, cf. \refde{Hhats}.

The content of the following theorems (\ref{theo_comparisongroups}, \ref{coro_ChhatVsArakelov}, \ref{theo_comparisonaction}, \ref{theo_comparisonpushforward}) can be paraphrased as follows: given a commutative diagram in some triangulated category,
$$\xymatrix{
B[-1] \ar[d]^{b[-1]} \ar[r] &
E \ar@{.>}[d]^e \ar[r] &
A \ar[d]^{a} \ar[r] &
B \ar[d]^{b} \\
B'[-1] \ar[r] &
E' \ar[r] &
A' \ar[r] &
B',
}$$
the map $e$ (whose existence is granted by the axioms of a triangulated category) is in general not unique. The unicity of $e$ is guaranteed if the map
\eqn \mylabel{eqn_EAB}
\Hom (E, A'[-1]) \r \Hom (E, B'[-1])
\xeqn
is onto. In our situation, we are aiming at a canonical comparison between, say, the groups $\Hhat^*$ and $\Khat_*$. Both theories arise from distinguished triangles where two of the three vertices are the same, namely the one responsible for $K$-theory and the one for Deligne cohomology. Moreover, the map between them considered \emph{up to homotopy}, i.e., in the triangulated category $\SH$, is the Chern class that is independent of choices---as opposed to the Chern form, which does depend on the choice of a hermitian metric on the vector bundle in question. As we shall see, the non-formal surjectivity of \refeq{EAB} is assured by conditions \refit{conda} and \refit{condb} of \refth{comparisongroups} (or condition \refit{condaQ} if one neglects torsion). Luckily, it only consists of an injectivity condition for the regulator on the base scheme $S$, not on all schemes $X \in \Sm/S$. This is one of the places where working with the objects representing the cohomology theories we are interested in is much more powerful than working with the individual cohomology groups.

\theo \mylabel{theo_comparisongroups}
Let $S$ be a regular scheme over an arithmetic ring. We suppose that
\begin{enumerate}[(a)]
\item \mylabel{item_conda} $\ch: K_0(S) \r \HD^0(S) = \oplus_p \HD^{2p}(S, p)$ is injective, and
\item \mylabel{item_condb} $K_1(S)$ is the direct sum of a finite and a divisible group.
\end{enumerate}
For example, these conditions are satisfied for $S = B = \Z$, $\RR$, or $\CC$.
Then the following holds:
\begin{enumerate}[(i)]
\item \mylabel{item_pfcomparisongroupsi} Given any maps $s$, $d$ in $\Ho_\bullet(S)$ such that the right square commutes, there is a unique $\widehat s \in \End_{\Ho(S)}(\Shat)$ making the diagram commute:
$$\xymatrix{
\Deligne_s = \Omega_s \Deligne_s[-1] \ar[d]^{\Omega_s d} \ar[r] & \Shat \ar@{.>}[d]^{\widehat s} \ar[r] & S_* \ar[r]^(.4){\ch_S} \ar[d]^{s} & \Deligne_s[-1] \ar[d]^d  \\
\Deligne_s = \Omega_s \Deligne_s[-1] \ar[r] & \Shat \ar[r] & S_* \ar[r]^(.4){\ch_S} &  \Deligne_s[-1]
}$$
\item \mylabel{item_pfcomparisongroupsii} Likewise, given any $b$ and $d$ making the right half commute in $\SH(S)$, there is a unique $\widehat b \in \End_{\SH(S)}(\BGLhat)$ making everything commute:
$$\xymatrix{
\oplus_p \HD\twi p[-1] \ar[d]^{d[-1]} \ar[r] & \BGLhat \ar@{.>}[d]^{\widehat b} \ar[r] & \BGL \ar[r]^{\ch} \ar[d]^b & \oplus_p \HD \twi p \ar[d]^{d}  \\
\oplus_p \HD\twi p[-1] \ar[r] & \BGLhat \ar[r] & \BGL \ar[r]^{\ch} & \oplus_p \HD \twi p.
}$$
\item \mylabel{item_pfcomparisongroupsiii}
The afore-mentioned unicity results give rise to a canonical isomorphism, functorial with respect to pullback
\eqn \mylabel{eqn_comparisonKhatHhat}
\Khat_n(X) \cong \Hhat^{-n}(X / S)
\xeqn
for any $X \in \Sm / S$, $n \geq 0$. (The definition of $\Khat_n(X)$ in \lcs is only done for $X/B$ proper, but can be generalized to non-proper varieties using differential forms with logarithmic poles at infinity, as in \refde{Burgos}.) %Alternatively, the reader may take the right hand side of \refeq{comparisonTakeda} as a definition of arithmetic $K$-groups.)
\end{enumerate}
%\item \mylabel{item_pfcomparisongroupsiv}
Instead of \refit{conda} and \refit{condb}, let us suppose that
\begin{enumerate}[(a)]
\setcounter{enumi}{2}
\item \mylabel{item_condaQ} $\ch: K_0(S)_\Q \r \HD^0(S) = \oplus_p \HD^{2p}(S, p)$ is injective. For example, this applies to arithmetic fields and open subschemes of $\SpecOF$ for a number ring $\OF$.
\end{enumerate}
Then there is a canonical isomorphism
\eqn  \mylabel{eqn_comparisonKhatHhatQ}
\Khat_n(X)_\Q \cong \Hhat^{-n}(X / S)_\Q.
\xeqn
\xtheo

\pf
\refit{pfcomparisongroupsii}: Let us write $(-,-) := \Hom_{\SH(S)}(-,-)$ and $R := \oplus_{p \in \Z} \HD \twi p$. Then we have exact sequences
\eqn
\mylabel{eqn_Homgroups}
\xymatrix{
(R, R[-1]) \ar[r]^\alpha \ar[d] & (\BGL, R[-1]) \ar[d]^\beta \\
(R, \BGLhat) \ar[r] & (\BGL, \BGLhat) \ar[r] \ar[d] & (\BGLhat, \BGLhat) \ar[r]^{\delta} & (R[-1], \BGLhat) \\
&  (\BGL, \BGL) \ar[d]^\gamma \\
&  (\BGL, R).
}
\xeqn
We prove the injectivity of $\delta$ by showing that both $\alpha$ and $\beta$ are surjective. For any $\Omega$-spectrum $E \in \SH(S)$ whose levels $E_n$ are $H$-groups such that the transition maps \refeq{transitionGrass} induce surjections $\Hom_{\HoC}(\Gr_{d,r}, \Omega_s^m E_n) \r \Hom_{\HoC}(\Gr_{d',r'}, \Omega_s^m E_n)$ for $m=1$, $2$, $n \geq 0$, there is an exact sequence
$$0 \r \R^1 \varprojlim E^1_\Omega \r \Hom_{\SH}(\BGL, E) \r \varprojlim E^0_\Omega \r 0.$$
Here, for any group $A$, $A_\Omega$ is the projective system
$$A_\Omega : \ \ \dots A [[t]] \r A [[t]] \r A [[t]] \r \dots \r A [[t]],$$
with transition maps $f \mapsto (1+t) df/dt$ and $E^r := \Hom_{\SH}(S^r, E)$ for $r=0$, $1$ \cite[IV.48, 49]{Riou:Thesis}. This applies to $E = \BGL$ and $E = R$, cf. \refeq{surjectivityGrassmannian}:
$$
\xymatrix{
0 \ar[r] & \R^1 \varprojlim (K_1(S)_\Omega) \ar[d] \ar[r] & \End(\BGL) \ar[d]^\gamma \ar[r] & \varprojlim (K_0(S)_\Omega) \ar[d] \ar[r] & 0 \\
0 \ar[r] & \oplus_p \R^1 \varprojlim (\HD^{-1}(S)_\Omega)  \ar[r] & \Hom(\BGL, R)\ar[r] & \oplus_p \varprojlim (\HD^0(S)_\Omega) \ar[r] & 0.
}$$
The left hand upper term is $0$ by assumption \refit{condb} and the vanishing of $\R^1 \varprojlim A_\Omega$ for a finite or a divisible group $A$ \cite[IV.40, IV.58]{Riou:Thesis}. The right hand vertical map $\varprojlim \ch$ is injective by assumption \refit{conda} and the left-exactness of $\varprojlim$. Hence $\gamma$ is injective, so $\beta$ is onto.

The surjectivity of $\alpha$ does not make use of the assumptions \refit{conda}, \refit{condb}. Indeed,
$$\Hom(\BGL, R[-1]) = \prod_{q \in \Z} \Hom(\HB \twi q, R[-1]) \stackrel{\text{\ref{theo_Delignespectrum}\refit{proofii}}}= \prod_q \HD^{-1}(S).$$
Given some $x \in \HD^{-1}(S)$, pick any representative $\xi \in \ker (\Deligne_{1}(S) \r \Deligne_0(S))$ and define $y: \HD\twi q \r R$ to be the cup product with $\xi$. Then $\alpha(y) = x$.

\refit{pfcomparisongroupsi}: we need to establish the injectivity of the map in the first row:
\eqn \mylabel{eqn_stablevsunstable}
\xymatrix{
\End_{\Ho_\bullet(S)} (\Shat) \ar@{=}[d] \ar[r] & \Hom_{\Ho_\bullet(S)} (\Omega_s \Deligne_s[-1], \Shat) \ar@{=}[d] \\
%-----
\End_{\Ho_\bullet(S)} (\Omega_{\P}^\infty \BGLhat)  \ar@{=}[d] \ar[r] &  \Hom_{\Ho_\bullet(S)} (\Omega_{\P}^\infty \HD[-1], \Omega_{\P}^\infty \BGLhat)  \ar@{=}[d] \\
%-----
\Hom_{\SH(S)} (\Sigma_{\P}^\infty \Omega_{\P}^\infty \BGLhat, \BGLhat) \ar[d]_{\Sigma_{\P}^\infty \rightleftarrows \Omega_{\P}^\infty} \ar[r] &   \Hom_{\SH(S)} (\Sigma_{\P}^\infty \Omega_{\P}^\infty \HD[-1], \BGLhat) \ar[d] \\
\Hom_{\SH(S)} (\BGLhat, \BGLhat)  \ar@{>->}[r]^\delta &   \Hom_{\SH(S)} (\HD[-1], \BGLhat).
}
\xeqn
The counit map $\Sigma_{\P}^\infty \Omega_{\P}^\infty \r \id$ is an isomorphism when applied to $\BGL$, and $\HD$ (and thus $\HD[-1]$), since these two spectra are $\Omega$-spectra. Therefore, the same is true for $\BGLhat$. We are done by \refit{pfcomparisongroupsii}. % as follows from the exact sequence (for any $n \geq 0$) of $\Hom$-GROUPS ?? in $\Ho_\bullet(S)$
%$$\dots \r \Hom (S^n X_+, \Omega \Deligne_s [-1]) \r \Hom (S^n X_+, \Shat) \r \Hom (S^n X_+, \Z \x \Gr) \r \Hom (S^n X_+, \Deligne_s [-1])$$
%and its stable analogue.

\refit{pfcomparisongroupsiii}: we obtain the sought isomorphism as the following composition:
\eqnarra
\Hhat^{-n}(X / S) & := &  \Hom_{\SH(S)}(\Sigma_{\P}^\infty S^n \wedge X_+, \hofib (\BGL \stackrel{\id \wedge 1_{\cDz}}\lr \BGL \wedge \cDz)) \nonumber \\
%----
& = %\explainiso{\refit{pfcomparisongroupsii}}
& \Hom_{\SH(S)}(\Sigma_{\P}^\infty S^n \wedge X_+, \hofib (\BGL \stackrel{\ch}\lr  \oplus_{p}\cDz \twi p))
\mylabel{eqn_pf1}
\\
%----
& = & \Hom_{\HoC(S)}(S^n \wedge X_+, \hofib (\Z \x \Gr \stackrel{\ch_0}\lr \Deligne_s))
\mylabel{eqn_pf2}
\\
%----
& = %\explainiso{\refit{pfcomparisongroupsi}}
& \Hom_{\HoC(S)}(S^n \wedge X_+, \hofib (\Omega_s S_* \stackrel{\ch_S}\lr \Deligne_s))
%----
\mylabel{eqn_pf3}
\\
%----
& = & \Hom_{\HoC(S)}(S^n \wedge X_+, \hofib (\Omega_s S_* \stackrel{\ch_S} \lr \Deligne_s))
\nonumber
\\
& = & \Hom_{\HoCsect(S)}(S^{n+1} \wedge X_+, \hofib (S_* \r \Deligne_s[-1]))
%----
%\\
%& = & \pi_{n+1} \left (\hofib_{\Delta^\op \PSh_\bullet(\Sm / X)}(S_* \r \Deligne_s[-1])(X) \right)
%%----
\mylabel{eqn_pf4}
\\
& = & \pi_{n+1} \left (\hofib_{\Delta^\op \Sets_\bullet}(S_*(X) \stackrel{\ch_S}\r \Deligne_s[-1](X)) \right)
%----
\nonumber
\\
&  \explainiso{\refeq{comparisonTakeda}} & \Khat_n(X).
\nonumber
\xeqnarra

The canonical isomorphism \refeq{pf1} follows from \refit{pfcomparisongroupsii}: we can pick representatives of $\BGL$ and of $\ch : \BGL \r \oplus \HD \twi p$ (\refth{mapSH}\refit{pfmapSHii}) in the underlying model category $\Spt$. We will denote them by the same symbols. We get a diagram of maps in $\Spt := \Spt^{\P}(\Delta^\op \PSh_\bullet(\Sm/S))$
$$\xymatrix{
\hofib (\id \wedge 1_{\HD}) \ar[r] \ar@{.>}[d]^\alpha & \BGL \ar[r]^(.4){\id \wedge 1_{\HD}} \ar@{=}[d] & \BGL \wedge \cDz \ar@{.>}[d]^\ch \\
\hofib (\ch) \ar[r] & \BGL \ar[r]^{\ch} & \oplus_p \cDz \twi p.
}$$
The Chern character for motivic cohomology and \refth{Delignespectrum}\refit{proofiii} induce an isomorphism $\ch: \BGL \wedge \HD \cong \oplus_p \HD\twi p$ in $\SH(S)$. As $\SH(S)$ is triangulated, we get some (a priori non-unique) isomorphism $\alpha$ in $\SH(S)$. By \refit{pfcomparisongroupsii}, however, it is unique.

Similarly, the isomorphism \refeq{pf3} follows from \refit{pfcomparisongroupsi}: still using the above lift of $\ch$ to $\Spt$, $\ch_0 := \Omega^\infty_{\P} \ch$ is a map of simplicial presheaves. The isomorphism $\tau: \Z \x \Gr \cong \Omega_s S_*$ \refeq{isoGrWaldhausen} can be lifted to a map $\tilde \tau$ of presheaves
$$\xymatrix{
\hofib \ch_0 \ar[r] \ar@{.>}[d] & \Z \x \Gr \ar[r]^{\ch_0} \ar[d]_{\tilde \tau} & \Deligne_s \ar@{=}[d] \\
\hofib \ch_S \ar[r] & \Omega_s S_* \ar[r]^{\ch_S} & \Deligne_s
}$$
The right hand square may not commute in $\Delta^\op \PSh(\Sm/S)$, but it does in $\Ho_\bullet(S)$. By \refit{pfcomparisongroupsi}, the resulting isomorphism (in $\Ho_\bullet(S)$) between $\hofib_{\Delta^\op \PSh}(\ch_0)$ and $\hofib_{\Delta^\op \PSh}(\ch_S)$ is independent of the choice of $\tilde \tau$ and $\ch_0$.

In order to explain the canonical isomorphisms \refeq{pf2}, \refeq{pf4}, recall the following generalities on model categories: let
$$F: \category C \leftrightarrows \category D: G$$
be a Quillen adjunction and let a diagram
$\delta: d_1 \stackrel f \lr d_2 \leftarrow *$ in $\category D$ be given. The homotopy fiber of $f$ is a fibrant replacement of the homotopy pullback of $\delta$. If $\category C$ and $\category D$ are right proper and $d_1$ and $d_2$ are \emph{fibrant}, then the homotopy pullback agrees with the homotopy limit and $\holim G (\delta)$ is weakly equivalent to $G \holim (\delta)$. Finally, replacing any object in $\delta$ by a fibrant replacement yields a weakly equivalent homotopy fiber \cite[19.5.3, 19.4.5, 13.3.4]{Hirschhorn:Model}. Thus
\eqn \mylabel{eqn_hofibAdjunction}
\Hom_{\Ho(\category D)}(F(c), \hofib f) = \Hom_{\Ho(\category C)}(c, \hofib G(f)).
\xeqn
We apply this to the Quillen adjunctions
$$
\Delta^\op(\PSh_\bullet(\Sm / X))
\begin{array}{c}  \id \\ \leftrightarrows \\ \id \end{array}
\Delta^\op(\PSh_\bullet(\Sm / X))
\begin{array}{c}  \Omega_{\P}^\infty \\ \leftrightarrows \\ \Sigma^\infty_{\P} \end{array}
\Spt^{\P}(\PSh_\bullet(\Sm/X)).
$$
The leftmost category is endowed with the sectionwise model structure, then the Nisnevich-$\A$-local one, and the stable model structure at the right. These model structures %(and all intermediate ones used to define $\SH(S)$ %*in \refeq{SH})
are proper
\cite[II.9.6]{GoerssJardine}, \cite[3.2., p. 86]{MorelVoevodsky:A1}, %for the third [minor difference: they start with simplicial sheaves instead of presheaves], \cite[4.2]{Jardine:Motivic} for the fourth
\cite[4.15]{Jardine:Motivic}. %for the fifth).
The simplicial presheaf $\Deligne_s$ is fibrant with respect to the section-wise model structure, since it is a presheaf of simplicial abelian groups. Moreover, it is $\A$-invariant and has Nisnevich descent by \refth{DeligneCohomology}\refit{etaleDescent}. Therefore, it is fibrant with respect to the Nisnevich-$\A$-local model structure. Moreover, $\cDz$ is an $\Omega$-spectrum by \refle{DeligneOmega}, so it is a fibrant spectrum (any level-fibrant $\Omega$-spectrum is stably fibrant \cite[2.7]{Jardine:Motivic}).
For \refeq{pf2}, we may pick a fibrant representative of $\BGL$ (still denoted $\BGL$) such that $\Omega^\infty_{\P} \BGL =: V$ is weakly equivalent to $\Z \x \Gr$. Again using \refit{pfcomparisongroupsi}, the homotopy fibers of $\Omega^\infty_{\P} (\ch): V \r \Deligne_s$ and of $\ch_0 : \Z \x \Gr \r \Deligne_s$ are canonically weakly equivalent.
%Secondly, $\BGL$ is stably fibrant THIS IS WRONG (REFEREE: \Z \x \Gr does NOT REPRESENT K-theory in Nisnevich sheaves, only so after inverting $\A$ : its terms $\Z \x \Gr$ are simplicially constant, i.e., fibrant in the Nisnevich model structure (since any representable presheaf is a Nisnevich sheaf) and $\A$-invariant, since $K_*(X) \cong K'_*(X)$ for all $X \in \Sm/S$ by the regularity of $S$.  Moreover, $\BGL$ is an $\Omega$-spectrum by the projective bundle formula for $K$-theory.
Finally, the $S$-construction presheaf $S_*$, cf.\ \refeq{Sconstruction}, is $\A$-invariant (since $K_*(X) \cong K'_*(X)$ for all $X \in \Sm/S$ by the regularity of $S$) and Nisnevich local for all regular schemes \cite[Thm. 10.8]{ThomasonTrobaugh}, and consists of Kan simplicial sets by definition. Hence $S_*$ is a fibrant simplicial presheaf in the $\A$-model structure. Therefore, \refeq{pf2}, \refeq{pf4} are fibrant, so these isomorphisms follow from \refeq{hofibAdjunction}.

The statement with rational coefficients is similar: one replaces $S_*$, which is given by simplicial chains in the topological realization of the $S$-construction, by its version with rational coefficients. Likewise, one replaces $\BGL$ by its $\Q$-localization (using the additive structure of $\SH(S)$) %cf. e.g.\ Riou appendix A
$\BGL_\Q$. Then condition \refit{conda} gets replaced by \refit{condaQ} and \refit{condb} becomes unnecessary, since the groups $\R^1 \varprojlim A^\Omega$ encountered above vanishes for a divisible group $A$.
\xpf

\subsection{Adams operations}
\refth{comparisongroups} can colloquially be summarized by saying that any construction on $\Khat_*$ etc.\ that is both compatible with the classical constructions on $K$-theory and Deligne cohomology and canonical enough to be lifted to the category $\SH(S)$ (or $\HoC(S)$) is unique. We now use this to study Adams operations on arithmetic $K$-theory. In \refsect{action} below, this principle is used to identify the $\BGL$-module structure on $\BGLhat$.

The arithmetic $K$-groups are endowed with Adams operations
\eqn \mylabel{eqn_AdamsKhat}
\Psi^k_{\Khat} : \Khat_n (X)_\Q \r \Khat_n (X)_\Q.
\xeqn
This is due to Gillet and Soul\'e \cite[Section 7]{GilletSoule:CharacteristicII} for $n=0$ and to Feliu in general \cite[Theorem 4.3]{Feliu:Adams}. Writing
$$\Khat_n(X)^{(p)}_\Q := \{ x \in \Khat_n(X)_\Q, \Psi^k_{\Khat} (x) = k^p \cdot x \ \ \text{for all }k \geq 1\}$$
for the Adams eigenspaces, the obvious question
\eqn \mylabel{eqn_decompKhat}
\oplus_{p \geq 0} \Khat_n(X)^{(p)}_\Q \stackrel ? = \Khat_n(X)_\Q.
\xeqn
was answered positively for $n=0$ in \lc, but could not be solved for $n > 0$ by Feliu since the management of explicit homotopies between the chain maps representing the Adams operations becomes increasingly difficult for higher $K$-theory. In this section, we show that the above Adams operations agree with the natural ones on $\Hhat^*(X)_\Q$ and thereby settle the question \refeq{decompKhat} affirmatively.

Feliu establishes a commutative diagram of presheaves of abelian groups,
$$\xymatrix{
C_1 := N \widehat C_* \ar[r]^{\ch_1} \ar[d]^{\Psi^k} & \Deligne_* \ar[d]^{\Psi^k_\Deligne} \\
C_2 := \tilde \Z \widehat C_*^{\widetilde{\P[]}} \ar[r]^{\ch_2} & \Deligne_*.
}
$$
The Adams operation $\Psi^k_\Deligne$ is the canonical one on a graded vector space,
$$\Psi^k_\Deligne : \Deligne_* := \oplus_p \Deligne_*(p) \r \oplus_p \Deligne_*(p), \Psi^k = \oplus_{p} (k^p \cdot \id).$$
The complexes $C_i$ at the left hand side are certain complexes of abelian presheaves defined in \ocs They come with maps $\Omega_s S_* \r \cK(C_i)$ that induce isomorphisms $K_*\t \Q = \pi_*(\Omega_s S_*) \t \Q \r \H_*(C_i) \t \Q$, $i=1,2$. By means of these isomorphisms, $\Psi^k$ corresponds to the usual Adams operation on $K$-theory (tensored with $\Q$). Moreover, both maps $\ch_i$ induce the Beilinson regulator from $K$-theory to Deligne cohomology.

Recall also the definition of the arithmetic Chow group from \cite[Section 3.3]{GilletSoule:Arithmetic} in the proper case and \cite[Section 7]{Burgos:Arithmetic} in general. %alternatively: \cite[2.4]{GilletSoule:CharacteristicI}.
\mylabel{CHhatGS} In a nutshell, the group $\CHhatGS^p(X)$ is generated by arithmetic cycles $(Z, g)$, where $Z \subset X$ is a cycle of codimension $p$ and $g$ is a Green current for $Z$, i.e. a real current satisfying $\Fr_\infty^* g = (-1)^{p-1} g$ such that $\omega (Z, g) := -\frac 1 {2\pi i} \partial \overline \partial g + \delta_Z$ is the current associated to a $C^\infty$ differential form (and therefore an element of $\Deligne_0(p)(X)$). Here $\delta_Z$ is the Dirac current of $Z(\CC) \subset X(\CC)$.
In analogy to the relation of $\Khat^T_0(X)$ vs. $\Khat_0(X)$, we put
\eqn \mylabel{eqn_CHhat}
\CHhat^p(X) := \ker (\omega : \CHhatGS^p(X) \r \Deligne_0(p)(X)).\footnote{The group $\CHhat^p(X)$ is denoted $\CHhat^p(X)_0$ in \lcs}
\xeqn

\coro \mylabel{coro_ChhatVsArakelov}
Under the assumption of \refth{comparisongroups}\refit{condaQ}, the isomorphism $\Khat_n(X)_\Q \cong \Hhat^{-n}(X)_\Q$ is compatible with the Adams operations $\Psi^k_{\Khat}$ on the left and, using the Arakelov-Chern class established in \refth{hats}, the canonical ones on the graded vector space on $\Hhat^{-n}(X)_\Q \cong \oplus_{p \in \Z} \Hhat^{2p-n}(X, p)$. In particular, there are canonical isomorphisms
\eqnarra
\Khat_n(X)^{(p)}_\Q &= & \Hhat^{2p-n}(X, p), \mylabel{eqn_Khat1} \\
\CHhat^p(X)_\Q = \Khat_0(X)^{(p)}_\Q  & =  &\Hhat^{2p}(X, p), \mylabel{eqn_Khat2} \\
\oplus_{p \in \Z} \Khat_n(X)^{(p)}_\Q  & = & \Khat_n(X)_\Q. \mylabel{eqn_Khat3}
\xeqnarra
\xcoro
\pf
We write $\Omega_{s, \Q} A := \varinjlim C_* (\Omega |A|)$ for any pointed connected simplicial set $A$. Here, $|-|: \Delta^\op \Sets \rightleftarrows \Top : C_*$ is the usual Quillen adjunction, $\Omega$ is the (topological) loop space, the direct limit is indexed by $\Z^{> 0}$ ordered by divisibility, and the transition maps $\Omega |A| \r \Omega |A|$ are the maps that correspond to the multiplication in $\pi_1(A)$. Then $\pi_n \Omega_{s, \Q}(A) = (\pi_n \Omega_s (A)) \t_\Z \Q$ for all $n \geq 0$. The construction is functorial, so it applies to the simplicial presheaf $S_*$ and gives us a $\Q$-rational variant denoted $S_{*,\Q}$. The map $\Psi^k: C_1 \r C_2$ yields an endomorphism $\Psi^k_S \in \End_{\Ho(S)}(S_{*, \Q})$. Moreover, the maps $\ch_i$, $i=1,2$ mentioned above factor over $\ch_{i, \Q}: S_{*, \Q} \r \Deligne_s[-1]$ and the obvious diagram $\ch_1$, $\ch_2$, $\Psi^k_\Deligne$ and $\Psi^k_S$ commutes up to simplicial homotopy, i.e., in $\HoCsect(S)$, a fortiori in $\Ho(S)$. Therefore, by \ref{theo_comparisongroups}\refit{pfcomparisongroupsi} we obtain a unique map $\Psi^k_{\Shat} \in \End_{\Ho(S)}(\Shat_{*, \Q})$, where $\Shat_{*, \Q} := \hofib \ch_1 : S_{*, \Q} \r \Deligne_s[-1]$. By construction, both $\Psi^k_{\widehat S}$ and the canonical Adams structure maps $\Psi^k_\Deligne \in \End_{\Ho(S)}(\Omega_s \Deligne_s[-1])$ map to the same element in $\Hom_{\Ho(S)}(\Omega_s \Deligne_s[-1], (\widehat S_*)_\Q)$. On the other hand, looking at
$$\xymatrix{
\BGLhat_\Q \ar[r] \ar@{.>}[d]^{\Psi^k_{\BGLhat}} &
\BGL_\Q \ar[r] \ar[d]^{\Psi^k_{\BGL}} &
\BGL_\Q \wedge \HD  \ar[d]^{\Psi^k_{\BGL} \wedge \id} \ar[r]^\ch_\cong &
R := \oplus_p \HD  \twi p \ar[d]^{\Psi^k_{\Deligne}}
\\ %-----------------------
\BGLhat_\Q \ar[r] &
\BGL_\Q \ar[r] &
\BGL_\Q \wedge \HD  \ar[r]^\ch_\cong &
R
}$$
there is a unique $\Psi^k_{\BGLhat} \in \End_{\SH(S)_\Q}(\BGLhat_\Q) \stackrel \delta \rightarrowtail \Hom(R[-1], \BGLhat_\Q)$ that maps to the image of the canonical Adams operation on the graded object $R[-1]$. Using $\End_{\SH}(R[-1]) = \End_{\Ho}(\Omega \Deligne_s[-1])$ (compare the reasoning after \refeq{stablevsunstable}) we see that the Adams operations on $\BGLhat_\Q$ and on $\Shat_{*, \Q}$ agree which yields the compatibility statement using the definition of the comparison isomorphism \refeq{comparisonKhatHhatQ}. The isomorphism \refeq{Khat1} is then clear, as is \refeq{Khat3}, using \refeq{HBhatVsBGLhat}. \refeq{Khat2} is a restatement of \cite[Theorem 7.3.4]{GilletSoule:CharacteristicII}.
%
%$$\xymatrix{
%K_{n+1, \Q}^{(p)} \ar[r] \ar@{^{(}->}[d] &
%\HD^{2p-n-1}( p) \ar[r] \ar@{^{(}->}[d] &
%\Khat_{n, \Q}^{(p)} \ar[r] \ar@{^{(}->}[d] &
%K_{n, \Q}^{(p)} \ar[r] \ar@{^{(}->}[d] &
%\HD^{2p-n}( p) \ar@{^{(}->}[d] \\
%% -----------
%K_{n+1, \Q} \ar[r] \ar[d]^\cong &
%\oplus_p \HD^{2p-n-1}( p) \ar[r] \ar@{=}[d] &
%\Khat_{n, \Q} \ar[r] \ar[d]^\cong &
%K_{n, \Q} \ar[r] \ar[d]^\cong &
%\oplus_p \HD^{2p-n}( p) \ar@{=}[d] \\
%% -----------
%\H^{-n-1}_\Q \ar[r] \ar@{->>}[d] &
%\oplus_p \HD^{2p-n-1}( p) \ar[r] \ar@{->>}[d] &
%\Hhat^{-n}_\Q \ar[r] \ar@{->>}[d] &
%\H^{-n}_\Q \ar[r] \ar@{->>}[d] &
%\oplus_p \HD^{2p-n}( p) \ar@{->>}[d] \\
%% -----------
%% -----------
%\H^{-n-1, p} \ar[r] &
%\HD^{2p-n-1}( p) \ar[r]  &
%\Hhat^{-n, p} \ar[r] &
%\H^{-n, p} \ar[r] &
%\HD^{2p-n}(p).
%}$$
%The rows are exact and the composition of the three vertical maps in the non-central columns are isomorphisms, so the five-lemma shows \refeq{Khat1}. The isomorphism between the arithmetic Chow and $K$-group is due to Gillet and Soul\'e \cite[Theorem 7.3.4]{GilletSoule:CharacteristicII}. Finally, \refeq{Khat3} follows from \refeq{Khat1} and \refth{longexactArakelov}\refit{proof1_iv}.
\xpf

\subsection{The action of $K$-theory on $\widehat K$-theory}\mylabel{sect_action}

Recall from \refre{propertiesHBhat} that $\BGLhat$ is a $\BGL$-module, i.e., there is a natural $\BGL$-action
$$\mu: \BGL \wedge \BGLhat \r \BGLhat.$$
For any smooth scheme $f: X / S$, this induces a map called the \emph{canoncial $\BGL$-action} on $\Hhat$-groups:
\eqnarr
\H^n(X) \x \Hhat^m(X) & = & \Hom_{\SH(S)} (X_+, \BGL[n]) \x \Hom (X_+, \BGLhat[m]) \\
& \r & \Hom (X_+ \wedge X_+, \BGL \wedge \BGLhat [n+m]) \\
& \stackrel {\Delta^* \circ \mu_*} \lr & \Hom (X_+, \BGLhat [n+m]) = \Hhat^{n+m}(X).
\xeqnarr
Here $\Delta: X_+ \r X_+ \wedge X_+ = (X \x X)_+$ is the diagonal map.

\theo \mylabel{theo_comparisonaction}
Let $S$ be a regular base scheme satisfying Condition \refit{condaQ} of \ref{theo_comparisongroups}. Then, at least up to torsion, the canonical comparison isomorphism $\Khat_n(X) \cong \Hhat^{-n}(X)$ is compatible with the canonical $\BGL$-action on the right hand side and the $K_*$-action
$$K_*(X) \x \Khat_*(X) \r \Khat_*(X)$$
induced by the product structure on $\Khat^T_*(X)$ established by Gillet and Soul\'e (for $\Khat_0$) \cite[Theorem 7.3.2]{GilletSoule:CharacteristicII} and Takeda (for higher $\Khat^T$-theory) \cite[Section 6]{Takeda:Higher} on the left hand side.

Similarly, the pairing
$$\CH^n(X) \x \CHhat^m(X) \r \CHhat^{n+m}(X)$$
induced by the ring structure on $\CHhatGS^*(X)$ agrees, after tensoring with $\Q$, with the canonical pairing
$$\H^{2n}(X, n) \x \Hhat^{2m}(X, m) \r \Hhat^{2(n+m)}(X, n+m).$$
\xtheo

\pf
Before proving the theorem proper, we sketch the definition of the product on $\Khat_*^T$: instead of the $S$-construction, Takeda uses the Gillet-Grayson $G$-construction $G_*(-) := G_*(\Phat(-))$ of the exact category of hermitian vector bundles on a scheme (see p. \pageref{hatP}). There is a natural weak equivalence $G_*(T) \r \Omega_s S_*(T)$. In particular, $\pi_n(G_*(T)) = K_n(T)$ for any scheme $T$ and $n \geq 0$. This gives rise to a
canonical isomorphism
$$\Khat_n (X) = \pi_n \hofib_{\Delta^\op(\Sets)} (G_*(X) \stackrel {\ch_G}\lr \Deligne_s(X)).$$
(cf. \cite[Theorem 6.2]{Takeda:Higher}). %This is stated in \lcs for $n \geq 1$ for the bigger groups $\Khat^T_n(X)$, but the proof carries over to $\Khat_0(X)$, as well.
The advantage of the $G$-construction is the existence of a bisimplicial version $G^{(2)}_*$ of $G$-theory together with a weak equivalence $R: G_* \r G^{(2)}_*$ and a map  $\mu_G: G_*(X) \wedge G_*(X) \r G^{(2)}_*(X)$, so that the induced map $\pi_n(G_*(X)) \x \pi_m(G_*(X)) \r \pi_{n+m}(G_*(X))$ is the usual product on $K$-theory. Moreover, $\ch_G$ factors over $R$.

Consider the following diagram, where $\mu_\Deligne : \Deligne_s \wedge \Deligne_s \r \Deligne_s$ is the product (cf. \refsect{Deligne}) and the terms in the second column denote the homotopy fibers (with respect to the section-wise model structure) of the respective right-most horizontal maps:
$$\xymatrix{
\Omega_s(G \wedge \Deligne_s) \ar[r] \ar[d]^{\Omega_s \mu_\Deligne \circ \ch_G} &
G \wedge \Ghat \ar[r] \ar@{.>}[d] &
G \wedge G \ar[r]^{\id \wedge \ch_G} \ar[d]^{\mu_G} &
G \wedge \Deligne_s \ar[d]^{\mu_\Deligne \circ \ch_G}
\\ %----------
\Omega_s \Deligne_s \ar[r] \ar@{=}[d] &
\widehat {G^{(2)}} \ar[r] &
G^{(2)} \ar[r] &
\Deligne_s
\\ %-------
\Omega_s \Deligne_s \ar[r] &
\widehat {G} \ar[r] \ar[u] &
G \ar[u]^R \ar[r]^{\ch_G} &
\Deligne_s . \ar@{=}[u]
}$$
The lower right square is commutative (on the nose) according to \lcs The upper right square is commutative up to (a certain) homotopy \cite[Theorem 5.2]{Takeda:Higher}, so there is some dotted map such that the left-upper square commutes up to homotopy. This yields a map $\phi: G \wedge \Ghat  \r \Ghat $ in $\Ho_\bullet(S)$ fitting into the following diagram (in $\HoC(S)$):
\eqn \mylabel{eqn_unicityproductG}
\xymatrix{
G \wedge \Omega_s \Deligne_s \ar[r] \ar[d]^{\mu_\Deligne \circ \ch_G} &
G \wedge \Ghat  \ar[r] \ar@{.>}[d]^\phi &
G \wedge G \ar[r] \ar[d]_{\mu_G} &
G \wedge \Deligne_s \ar[d]^{\mu_\Deligne \circ \ch_G}
\\ %---------
\Omega_s \Deligne_s \ar[r] &
\Ghat  \ar[r] &
G \ar[r] &
\Deligne_s.
}
\xeqn
%compatible with the $G$-action maps (in $\Ho_\bullet(S)$), $G \wedge G \r G$ and $G \wedge \Deligne_s \r \Deligne_s$.
The $K_*$-action on $\Khat_*$ is induced by $\phi$. Thus, to prove the theorem, it is sufficient to show that the diagram
$$\xymatrix{
\Omega_{\P}^\infty (\BGL \wedge \BGLhat) \ar[r]^(.7){\cong} \ar[d]^{\Omega_{\P}^\infty \mu} &
G \wedge \Ghat  \ar[d]^\phi \\
\Omega_{\P}^\infty (\BGLhat) \ar[r]^{\cong} &
\Ghat
}$$
is commutative in $\HoC(S)$. Here the horizontal isomorphisms are the ones from \refth{comparisongroups}. For this, it is sufficient to show that the dotted map in \refeq{unicityproductG} is unique (in $\Ho_\bullet(S)$). The latter statement looks very much like \ref{theo_comparisongroups}\refit{pfcomparisongroupsi}. Indeed, it can be shown in the same manner, as we now sketch: again, one first does the stable analogue, namely the unicity of a map $\BGL \wedge \BGLhat \r \BGLhat$ in $\SH(S)$ making the diagram analogous to \refeq{unicityproductG} commute. To do so, the sequences in \refeq{Homgroups} are altered by replacing $\Hom(?, *)$ by $\Hom(\BGL \wedge ?, *)$ everywhere. For any $E \in \DMBei(S)$, we have
$$\Hom_{\SH(S)_\Q}(\BGL \wedge ?, E) = \prod_{p \in \Z} \Hom_{\SH(S)_\Q}(\HB \twi p \wedge ?, E) = \prod_{p} \Hom_{\SH(S)_\Q}(?\twi p, E)$$
since $\DMBei (S) \subset \SH(S)_\Q$ is a full subcategory. This applies to $E = \HD$ and $E = \BGL_\Q = \oplus_p \HB \twi p$. Therefore, both the surjectivity of $\alpha$ and the injectivity of $\gamma$ in \refeq{Homgroups} carries over to the situation at hand.\footnote{I need to restrict to $\Q$-coefficients, since I do not know how to compute $\BGL \wedge \BGL$.} Then, the unstable unicity statement mentioned above is deduced from the stable one.

The statement for the arithmetic Chow groups follows from this: $\CHhat^*(X)_\Q$ is a direct factor of $\Khat_0(X)_\Q$ in a way that is compatible with the action of the direct factor $\CH^*(X)_\Q \subset K_0(X)_\Q$, by the multiplicativity of the arithmetic Chern class $\Khat_0^T(X)_\Q \cong \oplus_p \CHhatGS^p(X)_\Q$ \cite[Theorem 7.3.2(ii)]{GilletSoule:CharacteristicII}. Similarly, the $\HB$-action on $\HBhat$ is a direct factor of the $\BGL_\Q$-action on $\BGLhat_\Q$.
\xpf

\subsection{Pushforward}
Let $f: X \r S$ be a smooth proper map. According to \refdele{pushforward},
$$\Hom(\BGL \r f_* f^* \BGL \stackrel{\tr^\BGL_f,\cong}\lr f_! f^! \BGL, \ \BGLhat)$$
defines a functorial pushforward
$$f_* : \Hhat^n(X) \r \Hhat^n(S)$$
and similarly
$$f_* : \Hhat^n(X, p) \r \Hhat^{n-2 \dim f}(S-\dim f),$$
where $\dim f := \dim X - \dim S$ is the relative dimension of $f$. We now compare this with the classical pushforward on arithmetic $K$ and Chow groups. Recall from \cite[Prop. 3.1.]{Roessler:Adams} the pushforward $f_* : \Khat_0^T(X) \r \Khat_0^T(S)$. This pushforward depends on an auxiliary choice of a metric on the relative tangent bundle. It should be emphasized that the difficulty in the construction of $f_*$ on the full groups $\Khat_0^T(X)$ is due to the presence of analytic torsion. We now show that its restriction to $\Khat_0(X)$ agrees with $f_* : \Hhat^0(X) \r \Hhat^0(S)$ in an important case. This shows that analytic torsion phenomena and the choice of metrics only concern the quotient $\Khat_0^T / \Khat_0$. See also \cite{BurgosFreixasLitcanu:Generalized} for similar independence results. 

%The pushforward of arithmetic $K$-theory groups $\hat K_n(-)$ defined by Roessler ($n = 0$) and Takeda ($n \geq 0$, \cite[Section 7.3.]{Takeda:Higher}) applies to smooth projective maps $f: X \r Y$ between arithmetic varieties (flat over $\Z$ and regular). A functoriality statement for analytic torsion, which constitute the technically most challenging part of this pushforward, has been given by Faltings and Ma \cite[Theorem 5.5.]{Faltings:Lectures}, \cite[(0.5)]{Ma:Formes}. These pushforwards The pushforward on arithmetic Chow groups \cite[Theorem 3.6.1]{GilletSoule:Arithmetic} is defined for all proper and smooth (over $\CC$) maps between arithmetic varieties. For the time being, no pushforward has been established for the higher arithmetic Chow groups \cite{BurgosFeliu:Higher}.

\theo \mylabel{theo_comparisonpushforward}
%Under the isomorphism , the pushforward on arithmetic Chow groups, tensored with $\Q$, agrees with the one on Arakelov motivic cohomology groups in the following two cases:

\begin{enumerate} [(i)]
 \item \mylabel{item_proofpushforwardFp} The pushforward $i_* : \Hhat^0(\Fp) = \H^0(\Fp) = \Z \r \Hhat^0(\Z) = \Z \oplus \RR$ is given by $(0, \log p)$.
 \item \mylabel{item_proofpushforwardTopDim}
Let $\OF$ be a number ring and $S \subset \SpecOF$ an open subscheme and let $f: X \r S$ be smooth projective. For any $n \in \Z$, the following diagram is commutative, where the right vertical map is the pushforward on Gillet-Soul\'e's arithmetic Chow groups \cite[Theorem 3.6.1]{GilletSoule:Arithmetic}, and the middle map is its restriction:
$$\xymatrix{
\Hhat^{2 (\dim X+n)}(X, \dim X+n) \ar[d]^{f_*} \ar[r]^{\cong}_{\ref{theo_comparisongroups}} & \CHhat^{\dim X+n}(X)_\Q \ar[d]^{f_*} \ar@{^{(}->}[r] & \CHhatGS^{\dim X+n}(X) \ar[d]^{f_*}
\\
\Hhat^{2+2n}(S, n+1) \ar[r]^{\cong, \ref{theo_comparisongroups}} & \CHhat^{n+1}(S)_\Q \ar@{^{(}->}[r] & \CHhatGS^{n+1}(S).
}$$
\item \mylabel{item_proofpushforwardTopDim2} Under the same assumptions, the following diagram commutes, where the right vertical map is the pushforward mentioned above and the middle one is its restriction. In particular, the restriction of the $\Khat^T_0$-theoretic pushforward to the subgroups $\Khat_0$ does not depend on the choice of the metric on the tangent bundle $T_f$ used in its definition.
$$\xymatrix{
\Hhat^{0}(X)_\Q \ar[d]^{f_*} \ar[r]^{\cong}_{\ref{theo_comparisongroups}} & \Khat_0(X)_\Q \ar[d]^{f_*} \ar@{^{(}->}[r] & \Khat^T_0(X)_\Q \ar[d]^{f_*} \\
\Hhat^{0}(S) \ar[r]^{\cong, \ref{theo_comparisongroups}} & \Khat_0(S)_\Q \ar@{^{(}->}[r] & \Khat^T_0(S).
}$$
%the pushforward $f_*: \Hhat^0(X) \r \Hhat^0(S)$ agrees with the restriction of the pushforward on $\Khat_0(-)$ to $\Khat^T_0(-)$ \cite[p. 671]{Takeda:Higher}, at least up to torsion. In particular, the restriction of the $\Khat^T_0$-theoretic pushforward to the subgroups $\Khat_0$ does not depend on the choice of the metric on the tangent bundle $T_f$ used in its definition.
%pushforward $f_*:  \Hhat^{2 (\dim X+n)}(X, \dim X+n) \r \Hhat^{2+2n}(S, n+1)$  agrees (via the isomorphism of \refth{comparisongroups}) with the map $f_* : \CHhat^{\dim X+n}(X)_\Q \r \CHhat^{n+1}(S)_\Q$ obtained by restricting the pushforward on arithmetic Chow groups $\CHhatGS$.
%\item \mylabel{item_proofpushforwardTopDim2} Under the same assumptions, the pushforward $f_*: \Hhat^0(X) \r \Hhat^0(S)$ agrees with the restriction of the pushforward on $\Khat_0(-)$ to $\Khat^T_0(-)$ \cite[p. 671]{Takeda:Higher}, at least up to torsion.
\end{enumerate}
\xtheo

In order to prove \refit{proofpushforwardTopDim}, we need some facts pertaining to the Betti realisation due to Ayoub \cite{Ayoub:Note}: for any smooth scheme $B / \CC$, let
$$-^\An: \Sm / B \r \AnSm / B^\An$$
be the functor which maps a smooth (algebraic) variety over $B$ to the associated smooth analytic space (seen as a space over the analytic space attached to $B$), equipped with its usual topology. (This functor was denoted $-(\CC)$ above.) The adjunction
$$\An^*: \PSh(\Sm/B, \CC) \leftrightarrows \PSh(\AnSm/B^\An, \CC): \An_*$$
between the category of presheaves of complexes of $\CC$-vector spaces on $\Sm/B$ and the similar category of presheaves on smooth analytic spaces over $B^\An$ carries over to an adjunction of stable homotopy categories:
\eqn \mylabel{eqn_adjunctionSHAn}
\An^*: \SH(B, \CC) \leftrightarrows \SH^\An(B^\An, \CC): \An_*.
\xeqn
%\Ho(\Spt^{\P_{B^\An}}(\PSh(\AnSm/B^\An, \CC))): \An_*.$$
%The right hand side is the homotopy category of $\P_{B^\An}$-spectra of presheaves with respect to the model structure obtained by replacing $\A$ by $\mathbb D^1$, the disk in $\CC^\An$, and the Nisnevich topology by the usual analytic topology
We refer to \cite[Section 2]{Ayoub:Note} for details and notation; we use $\P_{B^\An}$-spectra instead of $(\A_{B^\An} / {\Gm}_{B^\An})$-spectra, which does not make a difference.
Secondly, there is a natural equivalence
$$\phi_X: \SH^\An(X^\An, \CC) %= \Ho(\PSh(\text{Open}(X^\An), \CC))
\stackrel \cong \lr \D(\Shvv_\An(X^\An, \CC))$$
of the stable analytic homotopy category and
%involving the homotopy category of presheaves on the open subschemes of $B^\An$ with respect to the local model structure for the analytical topology, which is the same as
the derived category of sheaves (of $\CC$-vector spaces), for any smooth $B$-scheme $X$. Both this equivalence and \refeq{adjunctionSHAn} are compatible with the exceptional inverse image and direct image with compact support in the sense that
$${f^\An}^! \phi_S \An^* = \phi_X \An^* f^!, \ \ f^\An_! \phi_X \An^* = \phi_S \An^* f_!$$
for any smooth map $f: X \r S$ of smooth $B$-schemes \cite[Th. 3.4]{Ayoub:Note}. Here $f_!$ and $f^!$ are the usual functors on the stable homotopy category, %*\refeq{adjshriek}
while ${f^\An}^!$ and $f^\An_!$ are the classical ones on the derived category.

To show \refit{proofpushforwardFp}, we need the following auxiliary lemma. It is probably well-known, but we give a proof here for completeness.

\lemm \mylabel{lemm_exacttriang}
In a triangulated category, let $A \stackrel \alpha \r B \stackrel \beta \r C \stackrel \gamma \r A[1]$ and $A' \stackrel {\alpha'} \r B' \stackrel {\beta'} \r C' \stackrel {\gamma'} \r A'[1]$ be two distinguished triangles.
Consider the maps of $\Hom$-groups induced by $\alpha$, $\alpha'$ etc. We suppose that $\beta^*$ is onto and $\gamma^*$ is bijective, as shown:
$$\xymatrix{
\Hom (B, A') \ar[d]^{\alpha^*} \ar[dr]^{\alpha'_*} &
\Hom (C, B') \ar@{->>}[d]^{\beta^*} \ar[dr]^{\beta'_*} &
\Hom (A[1], C') \ar[d]^{\gamma^*, \cong} \ar[dr]^{\gamma'_*}
\\
\Hom(A, A') &
\Hom (B, B') &
\Hom (C, C') &
\Hom (A[1], A'[1]).
}$$
%
%$$\xymatrix{
%\Hom(A, A') &
%\Hom (B, A') \ar[l]^{\alpha^*} \ar[d]^{\alpha'_*} \\
%%-----
%&
%\Hom (B, B') &
%\Hom (C, B') \ar@{->>}[l] \ar[d]^{\beta^*} \\
%%----
%&
%&
%\Hom (C, C') &
%\Hom (A[1], C') \ar[l]^\cong \ar[d] \\
%%-----
%&
%&
%&
%\Hom (A[1], A'[1]).
%}$$
Then, for any $\xi \in \Hom(B, A')$, $(\alpha^* \xi)[1] = (\xi \circ \alpha)[1]$ agrees with the image of any lift of $\alpha'_* \xi$ in $\Hom(A[1], A'[1])$ under the above maps.
\xlemm
\pf
Consider the following diagram
$$\xymatrix{
B \ar[r]^\beta \ar[d]^\xi \ar@{}[dr]|{(1)}&
C \ar[r]^\gamma \ar@{.>}[d]^\upsilon \ar@{}[dr]|{(2)} &
A[1] \ar[r]^{\alpha[1]} \ar@{.>}[d]^{\zeta, \zeta'} \ar@{}[dr]|{(3)} &
B[1] \ar[d]^{\xi[1]}
\\ %------
A' \ar[r]^{\alpha'} &
B' \ar[r]^{\beta'} &
C' \ar[r]^{\gamma'} &
A[1].
}$$
By assumption, there is a map $\upsilon$ making the square $(1)$ commute. Next, there is a unique map $\zeta$ making the square $(2)$ commute. On the other hand, by the axioms of a triangulated category, there is a (a priori non-unique) map $\zeta'$ making both $(2)$ and $(3)$ commute. Therefore, $\zeta = \zeta'$. This implies the claim.
\xpf

%$$\xymatrix{
%\Hom(\one, \HBhat \twi 1) &
%\Hom (i_* i^* \one, \HBhat \twi 1) \ar[l]^\cong \ar[d]^\cong \\
%%-----
%&
%\Hom (i_* i^* \one, \HB \twi 1) &
%\Hom_U (\one[1], \HBei{U} \twi 1) \ar[l]^\cong \ar[d]^{!!!???\cong} \\
%%----
%&
%&
%\Hom_U (\one[1], \HD \twi 1) &
%\Hom_\Z (\one[1], \HD \twi 1) \ar[l]^\cong \ar[d]^\cong \\
%%-----
%&
%&
%&
%\Hom_\Z (\one[1], \HBhat (1)[3])
%}$$

\pf (of \refth{comparisonpushforward})
\refit{proofpushforwardFp}: let $i: \SpecFp \r S := \SpecZ \leftarrow U := \Spec \Z[1/p]: j$. Consider the triangles
$$S^0 \r i_* i^* S^0 \r j_! j^* S^0 [1] \r S^0 [1],$$
$$\BGLhat \r \BGL \stackrel \ch \r \oplus_p \HD \twi p \r \BGLhat[1].$$
The assumptions of \refle{exacttriang} are satisfied, as can be checked using \refeq{KhatTakeda}: the generator of $K_0(\Fp)$ lifts to $(p, \pm 1)$ under $K_1(U) = p^\Z \x \{ \pm 1 \} \twoheadrightarrow K_0(\Fp)$, which in turn gets mapped to $\log p \in \HD^1(\Q, 1) = \RR$ under the Beilinson (or Dirichlet) regulator, which agrees with the Chern class $\ch$ by \refth{comparisonRegulator}. Therefore, the pushforward
$i_* : \Hhat^0(\Fp) = \H^0(\Fp) = K_0(\Fp) = \Z \r \Hhat^0 (\Z) = \Khat_0(S) = \Z \oplus \RR$ is the map $(0, \log p)$, so it agrees with the classical $\Khat$-theoretic pushforward.

\refit{proofpushforwardTopDim}: put $d' := d+n$. We need to show the commutativity of the following diagram:
\eqn
\mylabel{eqn_pushforwardcompatible}
\xymatrix{
(\HB, f^! \HBhat\twi{n+1}) \ar@{=}[d] \ar[r]^{\widehat p}
&
(\HB, \HBhat\twi {d'}) \ar[r]^\cong
&
\CHhat^{d'}(X)_\Q \ar[dd]^{f_*}
\\ %---------------------------------------
(\HB, f_! f^! \HBhat\twi{n+1}) \ar[d]^{f_! f^! \r \id}
\\ %---------------------
(\HB, \HBhat\twi {n+1}) \ar[rr]^\cong
&
&
\CHhat^{n+1}(S)_\Q.
}
\xeqn
Here $\widehat p$ is the relative purity isomorphism $f^! \HBhat \twi 1 \cong f^* \HBhat \twi d$.

We may assume $n \geq 0$ since $\CHhat^{\leq 0}(S) = 0$. The group $\CH^{d'}(X)$ is finite for $n=0$ by class field theory \cite[Theorem 6.1]{KatoSaito:Global} and zero for $n > 0$. %it is finite for any scheme $V$ of dimension $d$ that is proper (of finite type) over $\Z$ and not annihilated by any prime $p$. Therefore $\CH^d(X)$ is finite,
Hence $\HD^{2d'-1}(X, d') \r \Khat_0(X)^{(d')}_\Q$ is onto, by \refth{longexactArakelov}. On the other hand, for dimension reasons, $\HD^{2d'-1}(X, d') = \HBetti^{2d'-2}(X, \RR(d'-1))$.
By definition, the pushforward in arithmetic Chow groups \cite[Thm. 3.6.1]{GilletSoule:Arithmetic} is compatible with
\eqnarra
f_*: %\HD^{2d-1}(X^\An, d-1) =
\HBetti^{2d'-2}(X^\An, \RR(d'-1)) & \r & \HBetti^{2n}(\CC^\An,\RR(n)) = \RR \mylabel{eqn_integrate}\\
\omega & \mapsto & \frac 1 {(2 \pi i)^{d-1}} \int_{X^\An} \omega. \nonumber
\xeqnarra
Let $C^*$ be the presheaf complex of real-valued $C^\infty$-differential forms on smooth analytic spaces. This is a flasque complex and its (presheaf) cohomology groups agree with Betti cohomology with real coefficients. The construction and properties of $\cDz$ (esp. \refth{DeligneCohomology}) carry over and yield a spectrum $\An_* (\cB)$ representing Betti cohomology. The maps of complexes of sheaves on the analytic site,
$$[\RR(p) \r \cO \r \Omega^1 \r \dots \r \Omega^{p-1}] \r \RR(p) \stackrel \sim \r C^*(p),$$
%(the second map is a quasi-isomorphism)
give rise to a map of spectra $\cDz(p) \r \An_* \cB(p)$. The rectangle \refeq{pushforwardcompatible} is functorial with respect to maps of the target spectrum. Thus, we can replace $\HBhat \twi {n+1}$ by $\An_* \cB(n+1)[2n+1]$ and $f_* : \CHhat^{d'}(X)_\Q \r \CHhat^{n+1}(X)_\Q$ by $f_* : \HBetti^{2d'-2}(X^\An, \RR(d'-1)) \r \HBetti^{2n}(\CC, \RR(n)) \stackrel {n=0} = \RR$. This settles our claim, since the adjointness map $f^\An_! {f^\An}^! \CC \r \CC$ does induce the integration map \refeq{integrate} \cite[Exercise III.20]{KashiwaraSchapira:Sheaves}.

\refit{proofpushforwardTopDim2} the diagram
$$\xymatrix{
K_1(X) \ar[r] \ar[d]^{f_*} & \HD^{-1}(X) \ar[d]^{f_* \circ (- \cup \Td T_f)} \ar[r] & \Khat_0(X) \ar[d]^{f_*} \ar[r] & K_0(X) \ar[d]^{f_*} \\
K_1(S) \ar[r] & \HD^{-1}(S) \ar[r] & \Khat_0(S) \ar[r] & K_0(X)
}$$
is commutative, see \cite[Section 7]{Takeda:Higher}. On the other hand, applying
$$\Hom_{\BGL-\Mod}(f_! f^* \BGL \stackrel{\tr^\BGL} \r f_! f^! \BGL \r \BGL, -)$$
to the triangle \refeq{keyproperty} yields a diagram which is the same, except that $K_*$ is replaced by $\H^{-*}$ and $\Khat_*$ by $\Hhat^{-*}$ (and their respective pushforwards established in \refdele{pushforward}). Indeed, the pushforward on Deligne cohomology induced by $\tr^\BGL$ (as opposed to $\tr^{\Beilinson}$) is the usual pushforward, modified by the Todd class. This is a consequence of \refth{RiemannRoch}.

Now, \refit{proofpushforwardTopDim2} is shown exactly as \refit{proofpushforwardTopDim}: the only non-trivial part is $\Khat_0(X)^{(d)}_\Q$, which is mapped onto by $\HD^{2d-1}(X, d)$, since $K_0(X)^{(d)}_\Q = \CH^d(X)_\Q = 0$.
\xpf

\rema
The same proof works more generally for $f_* : \Hhat^n(X, p) \r \Hhat^{n-2 \dim f}(S, p-\dim f)$, provided that $\H^n(X, p) = K_{2p-n}(X)^{(p)}_\Q \r \HD^n(X, p)$ is injective. For example, given a smooth projective complex variety $X$ of dimension $d$, a conjecture of Voisin \cite[11.23]{Voisin:HodgeII} generalizing Bloch's conjecture on surfaces satisfying $p_g = 0$ says that the cycle class map $K_0(X)^{(d-l)}_\Q \cong \CH^{d-l}(X)_\Q \r \HBetti^{2(d-l)}(X, \Q)$ is injective (or, equivalently, that the cycle class map to Deligne cohomology is injective) for $l \leq k$ if the terms in the Hodge decomposition, $\H^{p,q}(X)$ are zero for all $p \neq q$, $q \leq k$.
\xrema

\section{The Arakelov intersection pairing} \mylabel{sect_motivicduality}

Let $S = \Spec \Z[1/N]$ be an open, non-empty subscheme of $\SpecZ$, where $N = p_1 \cdot \hdots \cdot p_n$ is a product of distinct primes. We write $Log(N) := \sum_i \Z \cdot \log p_i \subset \RR$ for the subgroup ($\cong \Z^n$) spanned by the logarithms of the $p_i$.

In this section, we give a conceptual explanation of the height pairing by showing that it is the natural pairing between motivic homology and Arakelov motivic cohomology.

\subsection{Definition}
\defi
For $M \in \SH(S)$, put
\eqnarr
\H_0(M) & := & \Hom_{\SH(S)} (S^0, M) \\
\H_0(M, 0) & := & \Hom_{\SH(S)_\Q} (S^0, M_\Q).
\xeqnarr
The second group is called \emph{motivic homology} of $M$ (seen as an object of $\SH$ with rational coefficients): for $M \in \DMBei(S)$, $\H_0(M, 0) \cong \Hom_{\SH(S)_\Q}(\HB, M_\Q)$.
\xdefi

\defi \mylabel{defi_intersectionPairing}
Fix some $M \in \SH(S)$. The \emph{Arakelov intersection pairing} is either of the following two maps
\eqnarr
\H_0(M) \x \Hhat^0 (M) & \r & \Hhat^0 (S^0) = \Khat_0(S) = \Z \oplus \RR / Log(N), \\
\pi_M: \H_0(M, 0) \x \Hhat^2 (M, 1) & \r & \Hhat^2 (S^0, 1) = \Khat_0(S)^{(1)}_\Q = (\RR / Log(N)) \t \Q,\\
(\alpha, \beta) & \mapsto & \beta \circ \alpha.
\xeqnarr
%given by the composition of morphisms in $\SH(S)$ and $\SH(S)_\Q$, respectively.
\xdefi

\rema \mylabel{rema_propertiespairing}
\begin{enumerate}[(i)]
\item
The tensor structure on the category $\DMBei^c(S)$, the subcategory of compact objects of $\DMBei(S) \subset \SH(S)_\Q$, is rigid in the sense that the natural map $M \r M\dual{}\dual$ is an isomorphism for any $M \in \DMBei^c(S)$, where $M\dual := \IHom_{\DM_{\Beilinson}(S)} (M, \HB)$ \cite [15.2.4]{CisinskiDeglise:Triangulated}. This implies that the natural map $\Hom(M, N) \r \Hom(N\dual, M\dual)$ is an isomorphism for any two such motives. In particular $\H_0(M, 0) \cong \H^0(M\dual, 0)$, so the pairing can be rewritten as
\eqn
\mylabel{eqn_Arakperfect}
\H^0(M\dual, 0) \x \Hhat^2 (M, 1) \r \H^2(S, 1).
\xeqn
This is the shape familiar from other dualities, such as Artin-Verdier duality,
$$\H^0(\SpecZ, \mathcal F\dual) \x \H^3_c (\Spec \Z, \mathcal F(1)) \r \H^3(\Spec \Z, \mu_\ell) = \Q / \Z.$$
In this analogy, an \'etale constructible $\ell$-torsion sheaf $\mathcal F$ corresponds to a motive $M$ and \'etale cohomology with compact support gets replaced by Arakelov motivic cohomology. The pairing \refeq{Arakperfect} is conjecturally perfect when replacing $\HBhat$ by $\HBhatRR$, which is constructed in the same way, except that $\HB$ gets replaced by $\HBRR$, a spectrum representing motivic cohomology tensored with $\RR$. The implications of this conjecture and its relation to special $L$-values is the main topic of \cite{Scholbach:SpecialL}.
\item \mylabel{item_functoriality}
By definition, the intersection pairing is functorial: given a map $f: M \r M'$, the following diagram commutes:
$$\begin{array}{cccccl}
\pi_M: & \H^0(M, 0) & \x & \Hhat^{2}(M\dual, 1) & \lr & \RR \\
& \uparrow & & \downarrow & & \downarrow = \\
\pi_{M'}: & \H^0(M', 0) & \x & \Hhat^{2}({M'}\dual,1) & \lr & \RR.
\end{array}
$$
\end{enumerate}
\xrema

\begin{comment} \todo!!!
%-----------------------------
\lemm
This pairing has the following properties:
\begin{enumerate}[(i)]

\item \mylabel{item_duality} It is compatible with duality in the sense that the following diagram commutes:
$$\begin{array}{cccccl}
\pi^i_M: & \Hhat^i(M) & \x & \H^{-i}(DM)_\RR & \lr & \RR \\
& \downarrow & & \uparrow & & \downarrow = \\
\pi^{-i}_{DM}: & \H^i(M)_\RR & \x & \Hhat^{-i}(DM) & \lr & \RR \\
& \downarrow & & \uparrow & & \downarrow = \\
	& \Hw^i(M) & \x & \Hw^{-1-i}(DM) & \lr & \RR \end{array}
$$
Here the lower row pairing is the duality of weak Hodge cohomology (\todo{DETAILS}).
\end{enumerate}
\xlemm
\pf
\refit{duality}: remains to be proven
\xpf
%-----------------------------
\end{comment}

\subsection{Comparison with the height pairing}
For a regular, flat, and projective scheme $X / \Z$ of absolute dimension $d$, Gillet and Soul\'e have defined the \emph{height pairing} $\mu_{GS}$:
$$\xymatrix{
\CH^m(X)_0 \ar@{^{(}->}[d] \ar@{}[r]|{\x} & \CH^{d-m}(X)_0 \ar[r]^{\mu_{\Beilinson}} & \CHhat^1(S) \ar@{=}[d] \\
\CH^m(X) \ar@{}[r]|{\x} & \CHhat^{d-m}(X) \ar@{^{(}->}[d] \ar@{->>}[u] \ar[r]^\mu &\CHhat^1(S) \ar@{=}[d] \\
\CHhat^m_{GS}(X) \ar@{}[r]|{\x} \ar@{->>}[u] & \CHhat^{d-m}_{GS}(X) \ar[r]^{\mu_{GS}} &\CHhat^1(S)
}
$$
Here, $\CH^m(X)_0 := \ker \CH^m(X) \r \HD^{2m}(X, m)$ is the subgroup of the Chow group consisting of cycles that are homologically trivial at the infinite place. The pairing $\mu$ is uniquely determined by $\mu_{GS}$. It is given by
$$(Z, (Z', g')) \mapsto (Z \cdot Z', \delta_Z \wedge g')$$
where $Z$ and $Z'$ are cycles of codimension $m$ and $d-m$, $\delta_Z$ is the Dirac current, and $g'$ is a Green current satisfying the differential equation
$$\omega(Z', g') = -\frac 1 {2\pi i} \partial \overline \partial g' + \delta_{Z'} = 0.$$
See \cite[Sections 4.2, 4.3]{GilletSoule:Arithmetic}  for details. The pairing $\mu_{\Beilinson}$ is the height pairing defined by Beilinson \cite[4.0.2]{Beilinson:Height}. More precisely, Beilinson considered the group of homologically trivial cycles on $X \x_S \Q$, but we will focus on the case where the variety in question is given over the one-dimensional base $S$. %The pairing $\mu_{GS}$ of Gillet and Soul\'e factors over the one in the middle line.

We now give a very natural interpretation of the height pairing $\mu$ in terms of the Arakelov intersection pairing. Our statement applies to smooth schemes $X$, only, essentially because of the construction of the stable homotopy category, which is built out of presheaves on $\Sm/S$ (as opposed to regular schemes, say).

\theo
Let $S \subset \SpecZ$ be an open (non-empty) subscheme and let $f: X \r S$ be smooth and proper of absolute dimension $d$. For any $m$, let $n := m - \dim f = m-d+1$ and let $M = \M(X) \twi n = f_! f^! \HB \twi{n}$ be the motive of $X$ (twisted and shifted). Then the height pairing $\mu$ (tensored with $\Q$) mentioned above agrees with the Arakelov intersection pairing in the sense that the following diagram commutes:
$$\xymatrix{
\CH^m(X)_\Q \ar[d]_{\cong}^{\text{\ref{sect_motives}}} \ar@{}[r]|\x &
\CHhat^{d-m}(X)_\Q \ar[d]_{\cong}^{\text{\ref{coro_ChhatVsArakelov}}} \ar[r]^\mu &
\CHhat^1(S)_\Q \ar[d]^\cong \\
\H_0(M, 0) \ar@{}[r]|\x &
\Hhat^2(M, 1) \ar[r]^{\pi_M} &
\Hhat^2(S, 1).
}$$
\xtheo
\pf
We need to show that the following diagram is commutative. Here $\one := \HB$ is the Beilinson motivic cohomology spectrum, $\onehat := \HBhat$ is its Arakelov counterpart (\refde{widehatA}), and $(-, -)$ stands for $\Hom_{\DMBei(?)}(-,-)$, where the base scheme $?$ is $S$ or $X$, respectively. Every horizontal map is an isomorphism. The maps labelled $p$ and $\widehat p$ are relative purity isomorphisms $f^! \cong f^* \twi {d-1}$, applied to $\one$ and $\onehat$, respectively. The isomorphisms between the (arithmetic) Chow groups and (Arakelov) motivic cohomology are discussed in \refsect{motives} and \refcor{ChhatVsArakelov}.
\small
$$\xymatrix{
(\one, f_! f^! \one\twi {n}) \ar[r]^{p} \ar@{}[d] |{\x} &
(\one, \one\twi {m}) \ar@{=}[r] \ar@{}[d] |{\x} &
(\one, \one\twi {m}) \ar[r] \ar@{}[d] |{\x} &
% K_0(X)^{(d-1-m)}_\Q \ar@{}[d] |{\x} \ar[r]^{\ch}&
\CH^{m}(X)_\Q \ar@{}[d] |{\x}
\\ % ----------------------
(f_! f^! \one\twi {n}, \onehat\twi 1) \ar@{}[rddd] |{(1)} \ar[ddd]^{\pi_M} \ar[r]^{p} &
(\one\twi {m}, f^! \onehat\twi 1) \ar[r]^{\widehat p} \ar[d]^\circ \ar@{}[rd] |{(2)} &
(\one\twi {m}, \onehat\twi d) \ar[d]^\circ \ar[r] \ar@{}[rd] |{(3)}  &
%\Khat_0(X)^{(1+m)}_\Q \ar[d]^\mu \ar[r]^{\widehat \ch} \ar@{}[rd] |{(5)}  &
\CHhat^{d-m}(X)_\Q \ar[d]^\mu
\\ %---------------------------------------
&
(\one, f^! \onehat\twi{1}) \ar@{=}[d] \ar[r]^{\widehat p} \ar@{}[rrdd] |{(4)} &
(\one, \onehat\twi d) \ar[r] &
%\Khat_0(X)^{(d)}_\Q  \ar[r]^{\widehat \ch}&
\CHhat^{d}(X)_\Q \ar[dd]^{f_*}
\\ %---------------------------------------
&
(\one, f_! f^! \onehat\twi{1}) \ar[d]^{f_! f^! \r \id}
\\ %  ----------------
(\one, \onehat\twi 1) \ar@{=}[r] &
(\one, \onehat\twi{1}) \ar[rr] &
&
%\Khat_0(S)^{(1)}_\Q \ar[r]^{\widehat \ch} &
\CHhat^{1}(S)_\Q
}$$
\normalsize
The commutativity of $(1)$ is a routine exercise in adjoint functors. The commutativity of $(2)$ is obvious. %The commutativity of $(5)$ follows from the multiplicativity of Gillet-Soul\'e's arithmetic Chern character $\widehat \ch : \widehat K_0^{GS}(X)_\Q \cong \oplus_p \CHhat^p_{GS}(X)_\Q$ \cite[Theorem 7.3.2(ii)]{GilletSoule:CharacteristicII}.
The commutativity of $(3)$ and $(4)$ is settled in Theorems \ref{theo_comparisonaction} and \ref{theo_comparisonpushforward}.
\xpf

\exam
Using \refre{propertiespairing}\refit{functoriality}, we can also describe the baby example of the Arakelov intersection pairing for $M = \M(\Fp)$: according to \refth{comparisonpushforward}\refit{proofpushforwardFp}, it is given by
$$\xymatrix{
\H_0(\Fp) \ar@{}[r]|\x &
\Hhat^0(\Fp) = \Z \ar[d]^{i_*}_{(0, \log p)} \ar[r]^{\pi_{\Fp}} &
\Hhat^0(\Z) = \Z \oplus \RR \ar@{=}[d] \\
\H_0(\Z) = \Z \ar[u]_{i^*}^\cong \ar@{}[r]|\x &
\Hhat^0(\Z) \ar[r]^{\pi_\Z} &
\Hhat^0(\Z) = \Z \oplus \RR
}$$
Using \refth{comparisonaction}, the bottom row is the obvious multiplication map. %using that $\Z \subset \Q$. Also note that, $1 \in \Z$ is induced by the unit of the ring spectrum $\BGL$, so it clearly acts as the identity on $\Z \oplus \RR$.
Therefore, $\pi_{\Fp}$ is given by $(1, 1) \mapsto (0, \log p)$.
\xexam

\bibliography{bib}
\end{document}